\documentclass{article}

\usepackage{arxiv}

\usepackage[utf8]{inputenc} 
\usepackage[T1]{fontenc}    
\usepackage{hyperref}       
\usepackage{url}            
\usepackage{booktabs}       
\usepackage{amsfonts}       
\usepackage{amsmath}        
\usepackage{nicefrac}       
\usepackage{microtype}      
\usepackage{ulem}           
\usepackage{float}          
\usepackage{lipsum}		
\usepackage{graphicx}
\usepackage{natbib}
\usepackage{doi}

\title{INI-VPINN: A Variational Physics-Informed Neural Network with Implicit Neumann and Interface Handling for Multi-Material Domains with Geometric Singularities}

\author{ \href{https://orcid.org/0000-0002-8323-2290}{\includegraphics[scale=0.06]{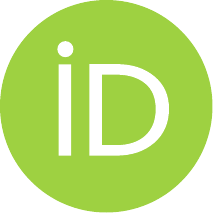}\hspace{1mm}Shayan Dodge}\thanks{Personal webpage: \url{https://shayandodge.github.io/}. Correspondence: \texttt{shayan.dodge@ing.unipi.it}.} \\
	DESTeC\\
	University of Pisa\\
	Pisa, Italy\\
	\texttt{shayan.dodge@ing.unipi.it} \\
	\And
	\href{https://orcid.org/0000-0002-7007-5759}{\includegraphics[scale=0.06]{orcid.pdf}\hspace{1mm}Alessandro Formisano} \\
	Department of Engineering\\
	University of Campania ``Luigi Vanvitelli''\\
	Aversa, Italy\\
	\texttt{Alessandro.formisano@unicampania.it} \\
    \And
	\href{https://orcid.org/0000-0003-1414-1114}{\includegraphics[scale=0.06]{orcid.pdf}\hspace{1mm}Sami Barmada} \\
	DESTeC\\
	University of Pisa\\
	Pisa, Italy\\
	\texttt{sami.barmada@unipi.it} \\
}

\renewcommand{\headeright}{Dodge et al.}
\renewcommand{\shorttitle}{INI-VPINN}

\hypersetup{
pdftitle={A template for the arxiv style},
pdfsubject={q-bio.NC, q-bio.QM},
pdfauthor={David S.~Hippocampus, Elias D.~Striatum},
pdfkeywords={First keyword, Second keyword, More},
}

\begin{document}
\maketitle

\begin{abstract}
We propose a new weak-form Physics-Informed Neural Network approach (named INI-VPINN). INI-VPINN naturally incorporates Neumann boundary and interface conditions into the variational formulation. It removes the need for additional loss terms or multiple subdomain networks. This framework employs compact support weighting functions and integration by parts to implicitly impose flux and continuity constraints. In this way, it implicitly ensures physical consistency across material boundaries.

The proposed method is tested on Poisson and Laplace problems with sharp interfaces and complex geometries. Results show that, compared with several other Physics Informed Neural Networks-based formulations, the INI-VPINN consistently achieves higher accuracy, smoother and faster convergence.

The proposed framework provides a general approach for solving multimaterial problems with complex geometries and mixed Neumann–Dirichlet boundary conditions using neural networks.

The implementation is publicly available in a GitHub repository (version v0.1.0).\footnote{\url{https://github.com/ShayanDodge/INI-VPINN}}
\end{abstract}


\section{Introduction}
\label{sec1}

Earlier developments in physics-guided neural network approaches for solving partial differential equations include the Deep Galerkin Method (DGM) \cite{sirignano2018dgm}, which laid important groundwork for subsequent methods. Building on these ideas, Physics-Informed Neural Networks (PINNs) were proposed for the resolution of Partial Differential Equations (PDEs) by Raissi et al. \cite{raissi2020hidden}. PINNs are a class of machine learning models that embed governing equations directly into the training process instead of relying purely on data. Equations are usually considered in the form of an additional loss function component at the output layer.

Thanks to several appealing features of this framework, it gained significant attention across a wide range of scientific and engineering fields, including:
fluid mechanics and Navier–Stokes flows~\cite{raissi2020hidden},  
heat transfer and thermal engineering~\cite{zhang2024physics},  
solid mechanics and inverse elastography~\cite{chen2023physics},  
geophysics and seismology~\cite{bandai2022forward},  
electromagnetics~\cite{barmada2025weak,dodge2025starpinn,barmada2024hybrid} and photonics~\cite{chen2022pinns_photonic_nanostructures},
cardiac electrophysiology and biomedicine~\cite{zhao2025pinns_physio_review},  
and even finance for option pricing problems~\cite{dhiman2023pinn_option_pricing}.

The key properties that contribute to the growing popularity of PINNs are: their \textit{mesh-free formulation}, which makes them well-suited for irregular domains~\cite{raissi2019physics}; \textit{physics-guided learning} that reduces, or completely removes, data requirements and improves extrapolation while maintaining physical consistency~\cite{karniadakis2021physics}; the \textit{lack of special discretization schemes}, because it integrates boundary and initial conditions into the loss function~\cite{lu2021deepxde};  their \textit{potential for high-dimensional problems}, mitigating the curse of dimensionality compared to classical solvers~\cite{sirignano2018dgm,yang2022deep}. 

Although PINNs have several notable strengths, they still face significant challenges, detailed as follows:

\begin{itemize}
    \item \textbf{Interface handling:} 
    In multi-material problems with sharp interfaces, standard PINNs struggle to satisfy the interface conditions. Many variants of the PINNs are introduced to mitigate interface issues like extended PINN (XPINN)~\cite{jagtap2020extended} with domain decomposition, interface PINNs (I-PINNs)~\cite{sarma2024ipinn} that combine domain decomposition with subdomain-specific activation functions, and conservative PINNs (cPINNs) ~\cite{jagtap2020conservative}. cPINNs still have limitations: they need one network per subdomain (thus increasing number of parameters and computational cost) and interface conditions are enforced via additional loss terms (making training sensitive to sampling density in the interface and to the associated loss weights).
    
    \item \textbf{Balancing multiple loss terms:} 
    In most of the problems in which PINNs are used, the need to enforce simultaneously PDE residuals, boundary (BC), initial (IC), and interface conditions arises. The multi-part loss often includes several terms with widely varying magnitudes, so that dominant ones frequently drown out the effect of smaller components in the total loss during the training~\cite{wang2022when}. To mitigate this imbalance in convergence of different loss terms many solutions have been introduced, e.g. \textit{a new gradient descent algorithm} using the eigenvalues of the Neural Tangent Kernel (NTK) to rescale and balance the convergence rates of different loss components~\cite{wang2022when} or a \textit{Self-adaptive point-wise weighting} learing a nonnegative weight per training point (IC/BC/residual) to balance the role of other terms in loss function~\cite{mcclenny2023self}. Although these strategies reduce sensitivity to hand-tuned loss weights, they often introduce more complexity and significant computational overhead. As a matter of fact, NTK-based methods require costly eigenvalue computations while self-adaptive schemes require learning many additional parameters. In addition, the risk of instability or bias is high since they do not guarantee global balance throughout training and dominant terms may still re-emerge as the solution evolves, leading to instabilities or slow convergence. Finally, such methods are often characterized by the correct training dynamics, but they do not fundamentally enhance the capacity of the neural network to resolve discontinuities or heterogeneous material discontinuities.

    \item \textbf{Complex geometries:} 
    Despite being mesh-free, PINNs still face challenges in irregular or multiscale domains, particularly those with sharp corners (singularities), where point sampling and BCs enforcement are non-trivial tasks. The most powerful approach to handle singularities is Variational PINNs with domain decomposition (hp-VPINN)~\cite{kharazmi2021hp}.
\end{itemize}

To address the above mentioned limitations of PINNs, several studies have proposed alternative approaches that impose the PDE in a weak form rather than a strong form. The Deep Energy Method (DEM)~\cite{Samaniego2020_energyML} introduced the idea of directly minimising the total energy function, thereby reducing numerical instabilities by avoiding high-order derivatives. Building on this perspective, Kharazmi et al.~\cite{KharazmiZhangKarniadakis2019_VPINN} proposed Variational Physics-Informed Neural Networks (VPINNs) as a Petrov-Galerkin extension of PINNs. VPINNs, unlike strong-form PINNs, enforce the weak formulation by projecting the residual onto global test functions and integrating over the entire domain. In parallel, Khodayi-Mehr et al.~\cite{KhodayiMehrZavlanos2020_VarNet} extended the idea in VarNet by employing locally supported, FE-like test functions for improved adaptivity; however, its accuracy depends on the chosen support size, and it lacks systematic refinement strategies. Other researchers have further generalized the weak enforcement principle; for instance, Weak Adversarial Networks (WANs)~\cite{ZangBaoYeZhou2020_WAN} cast the problem as a min–max game between a solution network and an adversarial test network. In this framework, the solution network attempts to minimize the weak form error by predicting a more precise solution. In contrast, the adversarial test network aims to maximize it by changing the test function. Finally, the Weak Form Theory-Guided Neural Network (TgNN-wf) ~\cite{XuZhangRongWang2021_TgNNwf} combined weak forms with theory-guided test functions and adaptive weight balancing.

 Later developments focused on enhancing efficiency and flexibility. The hp-Variational Physics-Informed Neural Network (hp-VPINN)~\cite{kharazmi2021hp} builds on the VPINN framework by combining variational loss enforcement with domain decomposition and bringing in the classical hp-refinement strategy from finite element methods to handle singularities and discontinuities.  Convolutional Variational PINN (cv-PINNs) ~\cite{LiuWu2023_cvPINN} accelerated this process by embedding variational operations into convolutional filters. 

Additional variants targeted specific challenges, such as: imposing Dirichlet boundary conditions using Approximate Distance Function (ADF) strategies and a weak imposition via Nitsche’s method~\cite{BerroneCanutoPintoreSukumar2023_DirichletBC_PINN_VPINN}; efficient handling of complex geometries in the Meshfree VPINN (MF-VPINN)~\cite{BerronePintore2024_MFVPINN} which employed locally supported patch-based test functions that avoided the need for a global mesh. More recently the Robust VPINNs (RVPINNs)~\cite{RojasMaczugaMunozMatutePardoPaszynski2024_RVPINNs} reformulated the loss function using the minimum-residual principle, making the method less sensitive to the choice of test functions. Overall, these advances demonstrate a clear evolution in PINN, with a significant focus on variational approaches that offer improved convergence behaviour, flexibility, and robustness.

Despite the substantial progress made by weak-form extensions of PINNs that alleviate many limitations of strong–form PINNs, they still leave two key issues largely unresolved: 
(i) the efficient enforcement of Neumann boundary conditions within the weak formulation, and 
(ii) the intrinsic enforcement of interface conditions in heterogeneous media.
These open issues motivated us to develop the Implicit Neumann and Interface VPINN (INI-VPINN), a framework that extends weak-form PINNs by naturally embedding Neumann boundary conditions and interface conditions into the variational residual. In this way, our approach avoids the need for ad-hoc penalty terms or separate subdomain networks.

\begin{table}[t]
\centering
\caption{Comparison of variational PINN frameworks.}
\renewcommand{\arraystretch}{1.15}
\setlength{\tabcolsep}{2.5pt}
\footnotesize
\begin{tabular}{p{2.8cm} p{5.5cm} p{4.7cm}}
\hline
\textbf{Method} & \textbf{Main Idea} & \textbf{Limitation} \\ \hline

\textbf{DEM}~\cite{Samaniego2020_energyML} & Minimizes total energy functional. & Only applicable when energy form exists.\\

\textbf{VPINN}~\cite{KharazmiZhangKarniadakis2019_VPINN} & Weak (Petrov–Galerkin) form with global weighting functions. & Global basis limits adaptivity in complex domains. \\

\textbf{VarNet}~\cite{KhodayiMehrZavlanos2020_VarNet} & Compact FE-like test functions. & Accuracy depends on patch size; no refinement rule. \\

\textbf{WAN}~\cite{ZangBaoYeZhou2020_WAN} & Min–max weak form with adversarial test network. & Costly optimization. \\

\textbf{hp-VPINN}~\cite{kharazmi2021hp} & Domain decomposition with $hp$-refinement for singularities. & Unable to implicitly handle Neumann and interface conditions. \\

\textbf{cv-PINN}~\cite{LiuWu2023_cvPINN} & Converts variational residuals to convolutions for GPU parallelism. & Less flexible for irregular geometries. \\

\textbf{RVPINN}~\cite{RojasMaczugaMunozMatutePardoPaszynski2024_RVPINNs} & Uses minimum-residual principle to find most efficient weighting function. & More expensive functional evaluation. \\

\textbf{INI-VPINN (ours)} & Implicit Neumann and interface enforcement via local integration; accurate for multi-material domains. & Slightly higher implementation complexity. \\

\hline
\end{tabular}
\end{table}

To ease reading, Table 1 reports a brief comparison of the different methods discussed here.

The main aim of this paper is then to demonstrate the ability of INI-VPINN in handling mixed boundary conditions — particularly Neumann ones — and non-homogeneous cases.
Section II introduces the INI-VPINN formulation and shows how Neumann and interface conditions are embedded in the variational residual. Section III describes the numerical implementation: compact-support test families, element-wise assembly. Section IV describes the neural network architecture and the training loss functions. Section V reports results on benchmarks (L- and T-shaped domains) and multi-material tests (nested square and circular inclusions), as well as a Poisson case with a compact source, with comparisons to vanilla PINN, VPINN with explicit Neumann enforcement, and cPINN. The paper concludes with a discussion of strengths, limitations, and future directions.

\section{Formulation of the INI-VPINN}
\label{sec2}

Let us consider an elliptic boundary-value problem described by the Poisson equation 
with spatially varying material properties on a domain 
$\Omega \subset \mathbb{R}^2$, with the boundary $\partial\Omega$ decomposed into Dirichlet 
and Neumann parts, $\partial\Omega = \Gamma_D \cup \Gamma_N$:

\begin{equation}
\left\{
\begin{aligned}
&-\nabla\cdot\!\big(\kappa(\mathbf{x}) \nabla u(\mathbf{x})\big) = f(\mathbf{x}), 
&& \forall \mathbf{x}\in \Omega, \\[6pt]
& u = g_D, 
&& \text{on } \Gamma_D, \\[6pt]
& \partial_n u = g_N, 
&& \text{on } \Gamma_N.
\end{aligned}
\right.
\label{eq:Poisson_strong}
\end{equation}

Here, $u(\mathbf{x})$ is an unknown potential function, $f(\mathbf{x})$ is the source term function, $\kappa(\mathbf{x})>0$ describes the space dependency of the material property, and $g_D$ and $g_N$ specify the Dirichlet and Neumann boundary conditions, respectively. The unknown function $u(\mathbf{x})$ must be sought  in a Sobolev space $W^{m,2}(\Omega)$, with $m\ge 2$.

Starting from \eqref{eq:Poisson_strong} we derive the weak formulation by multiplying the Poisson equation times a weighting function $v(\mathbf{x})$ and integrating over the domain $\Omega$. This results in:

\begin{equation}
\int_{\Omega} \big(-\nabla\!\cdot(\kappa(\mathbf{x}) \nabla u)\big)\, v(\mathbf{x}) \,\mathrm{d}\Omega
= \int_{\Omega} f(\mathbf{x})\, v \,\mathrm{d}\Omega.
\end{equation}

The selection of the specific weighting functions $v(\mathbf{x})$ and their support is crucial for the accuracy and computational efficiency of the method. This will be discussed later in the paper.
After applying Green’s identity on the whole domain and splitting the boundary integral according to 
$\partial\Omega = \Gamma_D \cup \Gamma_N$, we obtain

\begin{equation}
\int_{\Omega} \kappa(\mathbf{x}) \nabla u \cdot \nabla v \,\mathrm{d}\Omega
- \int_{\Gamma_D} \kappa(\mathbf{x}) \partial_n u \, v \,\mathrm{d}\Sigma
- \int_{\Gamma_N} \kappa(\mathbf{x}) \partial_n u \, v \,\mathrm{d}\Sigma
= \int_{\Omega} f\, v \,\mathrm{d}\Omega.
\label{eq:Poisson_weak_1}
\end{equation}

Equation \eqref{eq:Poisson_weak_1} allows to relax the condition on the weak derivative requirement, so $u(\mathbf{x})$ can be sought in the Sobolev space $W^{m-1,2}(\Omega)$.

As shown in equation \eqref{eq:Poisson_weak_1}, two boundary integrals appear: one over the Dirichlet boundary $\Gamma_D$ and another one over the Neumann boundary $\Gamma_N$. 

To avoid calculating normal derivatives of the potential function $(\partial_n u)$ on the Dirichlet boundary $ \Gamma_D$, the weighting functions are chosen such that
they belong to $W^{m-1,2}_0(\Omega)$ when their support includes $\Gamma_D$. On the contrary, when their support does include $\Gamma_N$ the WF can be chosen in $W^{m-1,2}(\Omega)$.
This eliminates the boundary integral associated with $\Gamma_D$, i.e.

\begin{equation}
\int_{\Gamma_D} \kappa(\mathbf{x}) \partial_n u \, v \,\mathrm{d}\Sigma =0
\label{eq:Dirichlet_weak}
\end{equation}
and $\partial_n u$ does not need to be evaluated on $\Gamma_D$. At the same time, the Neumann BCs can be directly imposed in the weak form: in order to do so, the weighting functions are chosen not to vanish on $\Gamma_N$. In this way

\begin{equation}
\int_{\Gamma_N} \kappa(\mathbf{x}) \partial_n u \, v \,\mathrm{d}\Sigma \neq 0
\label{eq:Neumann_weak}
\end{equation}
and the Neumann condition $\partial_n u = g_N$ be enforced "naturally" through the boundary integral.

Thus, we obtain

\begin{equation}
\left\{
\begin{aligned}
&\int_{\Omega} \kappa(\mathbf{x}) \nabla u \cdot \nabla v \,\mathrm{d}\Omega
- \int_{\Omega} f\, v \,\mathrm{d}\Omega
= \int_{\Gamma_N} \kappa(\mathbf{x}) g_N \, v \,\mathrm{d}\Sigma,
&& \forall v \in W^{m-1,2}(\Omega), \\[6pt]
&u = g_D, \;\; v = 0,
&& \text{on } \Gamma_D, \\[6pt]
&\partial_n u = g_N,\;\; v \neq 0,
&& \text{on } \Gamma_N.
\end{aligned}
\right.
\label{eq:weak-global}
\end{equation}

Note that the weak formulation \eqref{eq:weak-global} not only 
handles the Neumann BC implicitly, but it is also capable of 
accommodating material discontinuities without the explicit enforcement 
of interface conditions.

Many numerical methods for the solution of PDEs are based on domain discretization, and often the single component of this discretization is called \textit{element}; on the contrary  in methods based on the projection on basis functions (for instance the meshless methods~\cite{Patel2020Meshless}) no discretization of the domain is performed. In these cases the functions on which the solution is projected do not necessarily have compact support.

The weighting functions considered  for IVI-PINN indeed have compact support as it will be shown later. As a consequence defining a set of weighting functions that achieve full coverage of the domain is equivalent to discretize the geometrical domain into a set of elements (using the common nomenclature of many classical numerical methods), each of them coincident with the support of the weighting functions. Consequently, from now on we will use the word element for each $V_k$ ($k = 1, \dots, N$) in which the domain $\Omega$ is discretized; for each element, a corresponding test function (or a set of them, as it will be shown later) $v_k$ is defined. 

The weak form can then be written in an element-wise representation:

\begin{equation}
\sum_{k=1}^{N} \int_{V_k} \kappa(\mathbf{x}) \nabla u \cdot \nabla v_k \,\mathrm{d}\Omega
- \sum_{k=1}^{N} \int_{V_k} f\, v_k \,\mathrm{d}\Omega
=\sum_{h=1}^{N_\Gamma} \int_{\Gamma_N} \kappa(\mathbf{x}) g_N\, v_h \,\mathrm{d}\Sigma
\label{eq:weak-local}
\end{equation}

in which $N_\Gamma$ is the number of elements having a common side with the boundaries characterized by Neumann conditions $(\Gamma_N)$; obviously  $N_\Gamma < N$ and \eqref{eq:weak-local} must hold  $\forall v_k \in W^{m-1,2}(\Omega)$.
This localized formulation \eqref{eq:weak-local} improves the convergence behaviour and provides better resolution of local features such as material discontinuities, sharp gradients, or corner singularities, as it will be shown later.

\section{Numerical Implementation}
\label{sec3}

\subsection{Weighting functions families}
The selection of weighting functions in the VPINN framework is not unique and depends on the specific problem.
Early studies on VPINN and hp-VPINN formulations \cite{KharazmiZhangKarniadakis2019_VPINN}, \cite{kharazmi2021hp} highlighted the flexibility of the approach, allowing for polynomial, trigonometric, and locally supported functions.
Standard VPINNs \cite{KharazmiZhangKarniadakis2019_VPINN} typically use orthogonal polynomial bases, while hp-VPINNs \cite{kharazmi2021hp}
introduce element-wise weighting functions through domain decomposition. The aim in this case is enhancing performance in problems with localized features or singularities. Our previous work on weak-form PINNs \cite{barmada2025weak} further demonstrated that the choice of weighting functions significantly affects accuracy, convergence, and boundary-condition treatment. Different families, such as Gaussian and compactly supported functions, show distinct behaviors regarding localization and interaction with domain boundaries.
In the INI-VPINN framework, different families of compact support weighting functions are available. For the sake of exposition, we describe below the three different families that we employed in the treatment of the benchmark problems reported here. Each of them is characterized by an oscillatory behavior determined by its rank \( n \), higher ranks corresponding to higher spatial frequencies, which can be adjusted according to the complexity of the problem.

In the remaining part of this section, for the sake of clarity and illustration, the 1D counterparts of the weighting functions used in Sect. 5 are described. Weighting functions are easily generalized to higher dimensional spaces by tensor product, as described below.

\begin{figure} 
    \centering
    \includegraphics[width=1\linewidth]{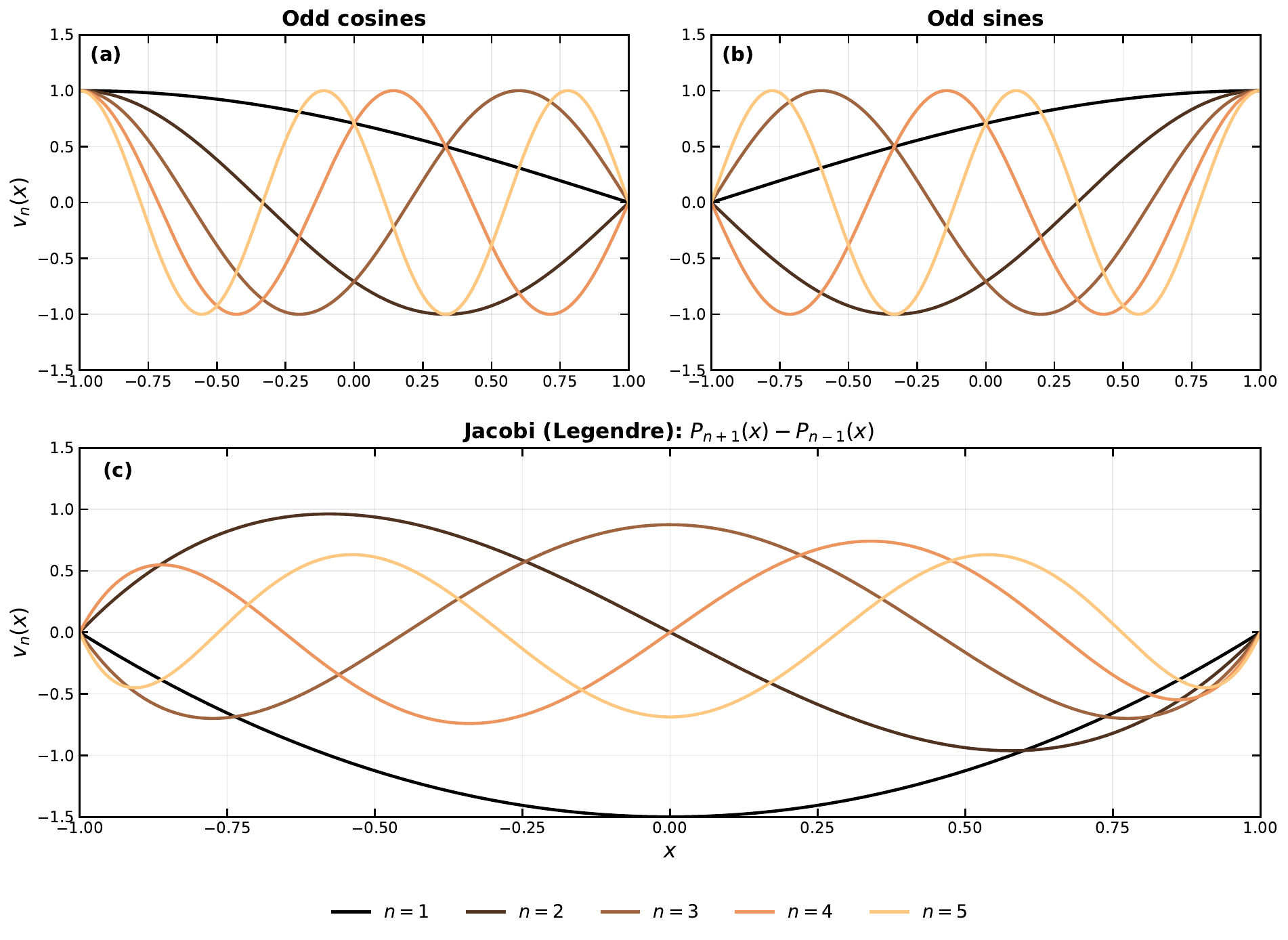} 
    \caption{Test function families: (a) odd cosines, (b) odd sines, and (c) Jacobi difference modes, shown for ranks $n= 1,\dots, 5$.}\label{fig:test_functions}
\end{figure}

\begin{itemize}
    \item \textbf{Odd cosine modes:}
For a general domain $x \in [a,b]$, the odd cosine weighting functions are defined by
\begin{equation}
v^{(\cos)}_n(x) 
= \cos\!\left( \frac{(2n-1)\pi}{2}\,\frac{x-a}{b-a} \right),
\qquad n=1,2,\dots
\end{equation}

These functions correspond to cosines whose wave number is an odd multiple of $\pi$; in Figure~\ref{fig:test_functions}(a) they are shown in the case of $[a,b]=[0,1]$. By construction, $v^{(\cos)}_n(a)=+1$ for all $n$ and $v^{(\cos)}_n(b)=0$ for all $n$, and the spatial frequency increases with the rank.

    \item \textbf{Odd sine modes:} 
For a general domain $x \in [a,b]$, the odd sine weighting functions are defined by
\begin{equation}
v^{(\sin)}_n(x) 
= (-1)^n \,\sin\!\left( \frac{(2n-1)\pi}{2}\,\frac{x-a}{b-a} \right),
\qquad n=1,2,\dots
\end{equation}

These functions correspond to sines whose wave number is an odd multiple of $\pi$; in Figure~\ref{fig:test_functions}(b) they are shown in the case of $[a,b]=[0,1]$. The factor $(-1)^n$ ensures that $v^{(\sin)}_n(b)=+1$ for all $n$, and $v^{(\sin)}_n(a)=0$ for all $n$, so the family shares a consistent endpoint normalization. They complement the cosine modes, forming an odd-harmonic trigonometric basis.

    \item \textbf{Jacobian modes (Difference of Legendre polynomials):} 
Let $P_n(x)$ denote the Legendre polynomial of degree $n$ defined on $[a,b]$. ~\cite{GuoShenWang2006_SpectralGalerkin}
The Jacobi-type difference functions are defined by

\begin{equation}
v^{(\mathrm{Jac})}_n(x) 
= P_{n+1}(x) - P_{n-1}(x),
\qquad n=1,2,\dots
\end{equation}

These functions are differences of consecutive Legendre polynomials. They retain orthogonality properties from the Jacobi family and 
produce polynomial-type oscillations whose complexity increases with $n$. 
The Jacobi functions are shown in Figure~\ref{fig:test_functions}(c). 
A key feature of these functions is that their values vanish at the domain boundaries, 
i.e., $v^{(\mathrm{Jac})}_n(a)=0$ and $v^{(\mathrm{Jac})}_n(b)=0$.
\end{itemize}

\subsection{Selection of the weighting functions}

Since the domain is divided into several non-overlapping elements, as explained before, the proper weighting functions must be chosen for each element among the three families introduced before; this choice is not trivial and is, on the contrary, fundamental for the accuracy of the solution.

In the INI-VPINN framework, a few key principles govern the choice of weighting functions, based on the boundary and interface conditions at each element. 

The first guideline is that for interior elements (i.e., those that do not touch the global boundaries), Jacobian functions are the best-suited option, since they vanish at both endpoints, making them consistent at interfaces, avoiding spurious boundary residuals, and providing stable convergence. 

Another guideline is for elements with a boundary in common with the domain boundary $\Gamma$: if the boundary condition is of the Neumann type, one should select between odd cosine and odd sine functions depending on the location of the Neumann boundary, choosing the family that does not vanish on that boundary section.

For Dirichlet boundaries, Jacobi functions or odd sine or cosines, depending on the point where the condition must be enforced, can be used to enforce vanishing value condition.
By setting the test functions $v$ equal to zero on Dirichlet boundaries, we avoid computing the boundary integrals in~\eqref{eq:weak-global} involving the normal derivative of $u$. However, this choice does not prescribe the actual value of $u$ on the boundary. The enforcement of $u=0$ (homogeneous case) or $u=g_d\neq 0$ (non-homogeneous case) must still be treated explicitly, e.g., via a lifting function or explicitly enforcing the boundary condition.

\begin{figure} 
    \centering
    \includegraphics[width=0.85\linewidth]{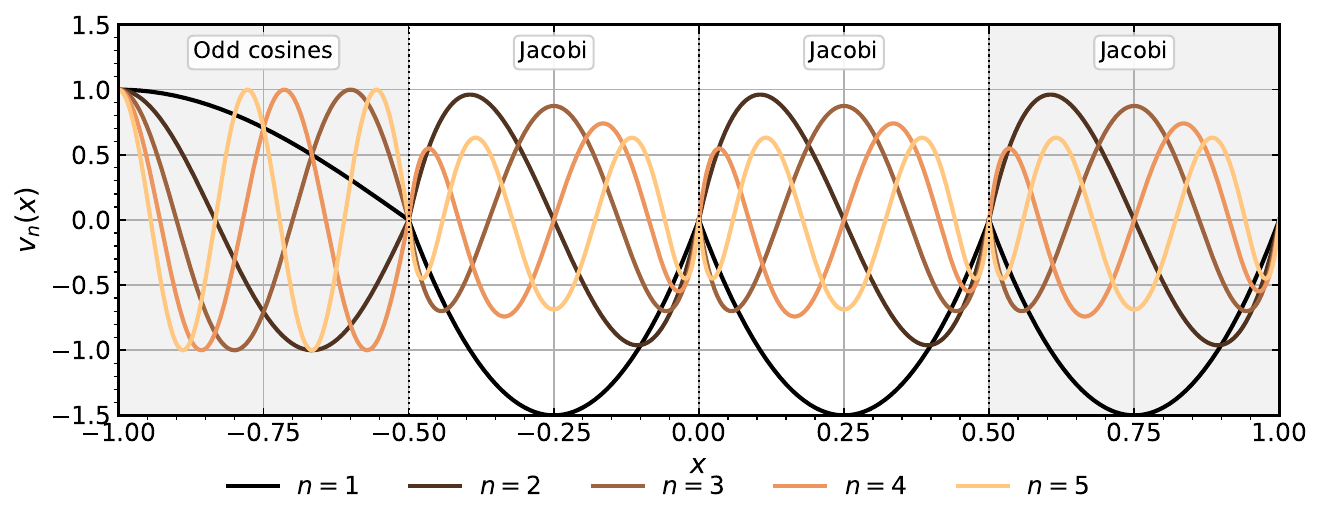} 
    \caption{Arrangement of test functions over a 1D domain: Odd cosine modes combined with Jacobi functions for Neumann-Dirichlet boundaries. Each family is shown for ranks $n=1$--$5$.}\label{fig:1d_tf}
\end{figure}

In order to show what explained before for a simple 1D configuration, in Figure~\ref{fig:1d_tf} the computational domain is divided into 
four non-overlapping elements; this is the case when one endpoint is characterized by Dirichlet type BC and the other one characterized by Neumann type BC.  As in 
in Figure~\ref{fig:1d_tf}, a combination of odd cosine modes 
(which satisfy $v^{(\cos)}_n(a)=+1$ at the Neumann side) and Jacobi modes 
is employed. 

Here, the truncation index $n$ denotes the number of basis functions retained in the approximation. In general $n$ is chosen large enough to achieve the desired accuracy. From the Petrov–Galerkin perspective, the truncation index $n$ plays the role of $p$ in the $ hp$-framework: increasing $n$ corresponds to $p$-refinement, while subdividing the domain corresponds to $h$-refinement.

\begin{figure}
\centering
  \includegraphics[width=.8\linewidth]{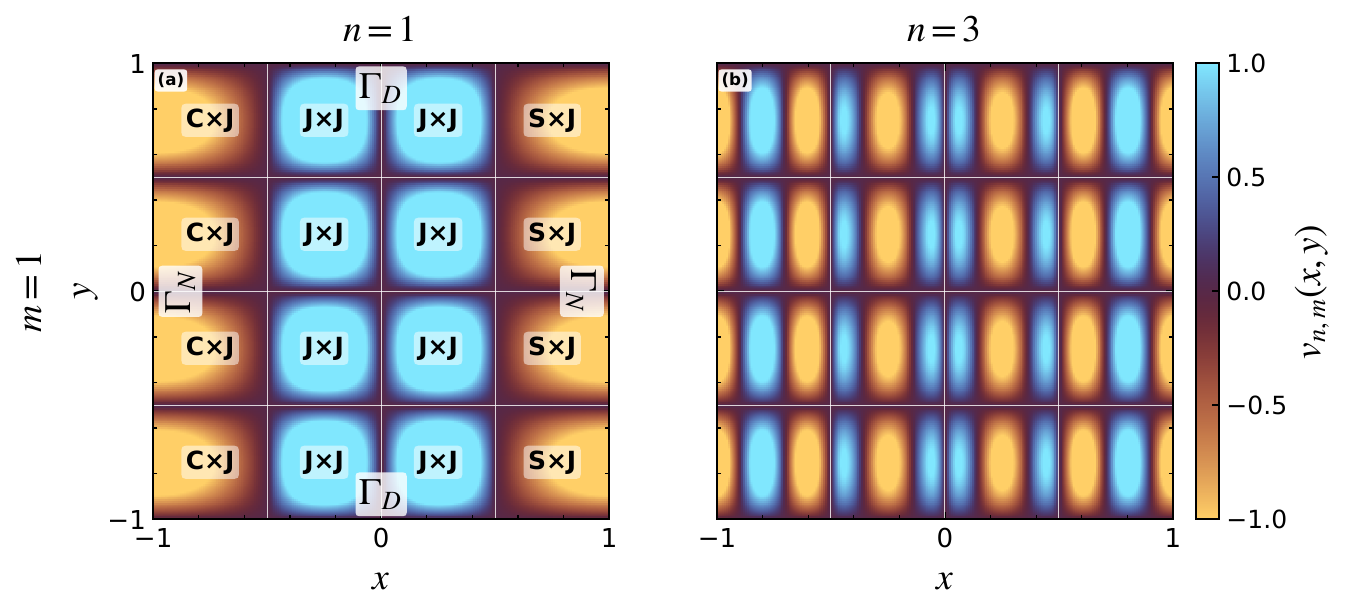}
\caption{Localized weighting function $v_{1,1}^{(k)}(x,y)$ 
    on a partitioned domain with mixed boundary conditions. These functions vanishes at the Dirichlet boundaries ($\Gamma_D$) and are normalised to one at the Neumann boundaries ($\Gamma_N$). Element labels indicate the chosen local basis type (J = Jacobi, C = odd cosine, S = odd sine). (a) $m=1, n=1$, (b}
\label{fig:2d_tf_square}
\end{figure}

In the case of a 2D problem, the domain is partitioned into $N$ smaller rectangular elements, indexed by 
$k = 1,\dots,N$, as illustrated in Figure~\ref{fig:2d_tf_square}. For each element, 
the type of weighting function (Jacobi, odd cosines, or odd sines) 
is selected based on its location. Independent choices are made for the 
$x$- and $y$-directions, with ranks $n=1,\dots,N^x$ and $m=1,\dots,N^y$, 
respectively. Consequently, each element is associated with a total of 
$N^x\times N^y$ tensor-product between suited 1D weighting functions. Let 
$v^{(x)}_{k,n}$ and $v^{(y)}_{k,m}$ denote the 1D test functions on element $k$ in the $x$- and $y$-directions, respectively. The corresponding tensor-product is defined as

\begin{equation}
v_{k,nm}^{(x,y)}(x,y) = v_{k,n}^{x}(x)\,v_{k,m}^{y}(y).
\label{eq:tensor_product}
\end{equation}

Also in the 2D case the truncation indexes $m,n$ denotes the number of basis functions retained in the approximation.

The INI-VPINN framework is not limited to the 2D problems; its extension to higher dimensions is achieved by extending the tensor product of the weighting functions \eqref{eq:tensor_product} by adding additional terms. However, the computational burden for the construction of the
tensor product and for the evaluation of the weak form on the 3D grid of points still needs further investigation.

\subsection{Implicit Enforcement of Interface Conditions}

\subsubsection{Sharp interface representation (Heaviside)}

Let us expand the Poisson equation as:

\begin{equation}
\nabla\!\cdot(\kappa\nabla u)=\kappa\,\nabla^2 u + \nabla \kappa \cdot\nabla u,
\label{eq:strong_kappa}
\end{equation}
and define $\Gamma_I$ as the interface between two different materials, respectively characterized by constants $\kappa_1$ and $\kappa_2$.

\begin{figure} 
    \centering
    \includegraphics[width=0.4\linewidth]{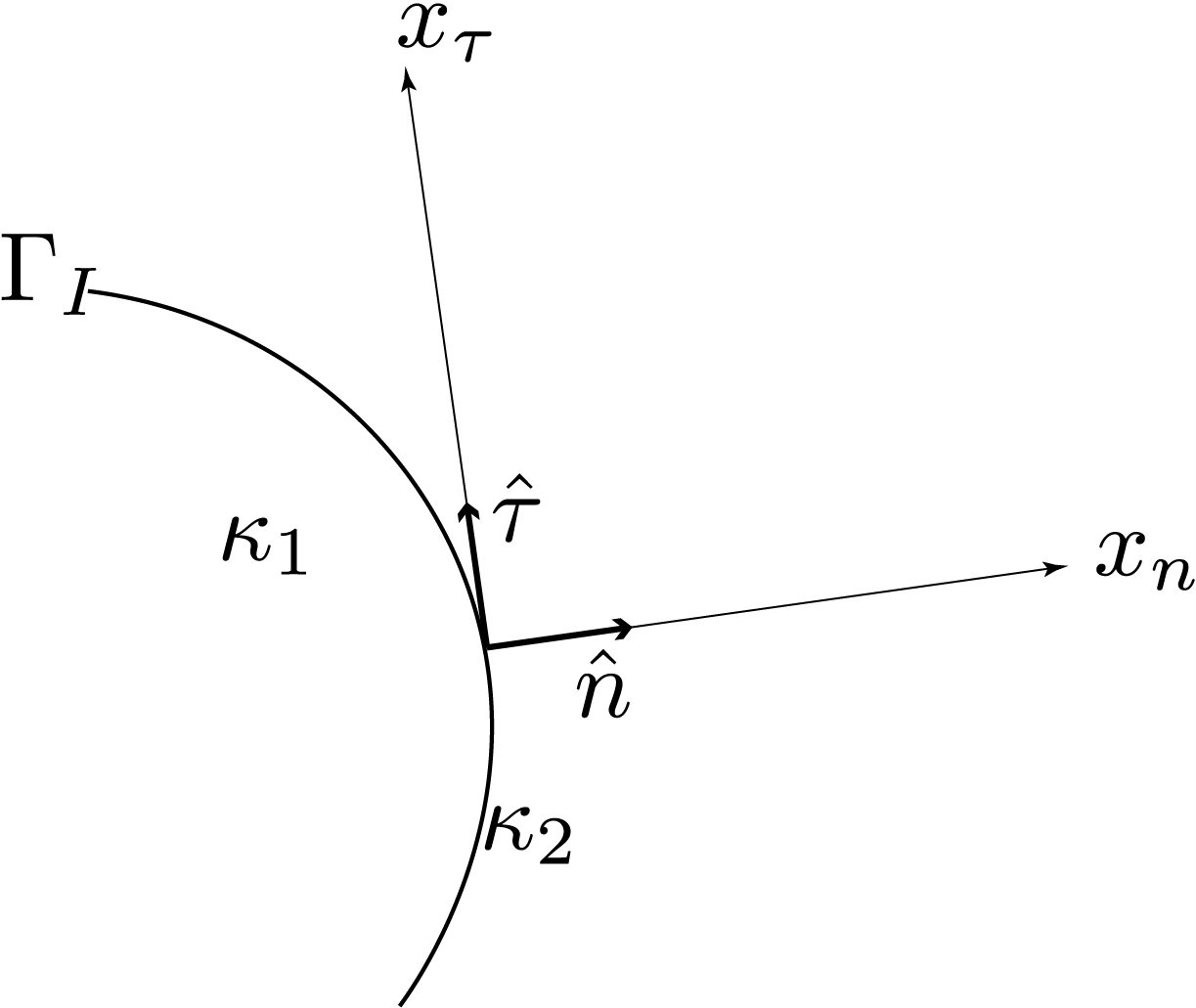} 
    \caption{Definition of the local coordinate system at the interface between two media.}
    \label{fig:interface_geometry}
\end{figure}

Figure \ref{fig:interface_geometry} shows the outline of an interface between two different materials: in each point $P\in \Gamma_I$ a local reference system $(x_n,x_\tau$) can be defined, whose axis are aligned with the normal versor $\hat n$ and the tangential versor $\hat \tau$ respectively (the extension to 3D geometries is trivial).

The most straightforward way to encode two different materials is to use the Heaviside step 
function. In terms of the coordinate systems described in Figure \ref{fig:interface_geometry} it reads as $H(x_n)$; in this way the change in the material property at the interface $\Gamma_I$ can be expressed as

\begin{equation}
\begin{aligned}
&\kappa(x_n)        = \kappa_1 + (\kappa_2-\kappa_1)\,H(x_n), \\[6pt]
&\nabla \kappa(x_n) = (\kappa_2-\kappa_1)\,\delta_D(x_n)\,\hat n,
\end{aligned}
\label{eq:H_delta}
\end{equation}

\begin{figure}
    \centering
    \includegraphics[width=0.85\linewidth]{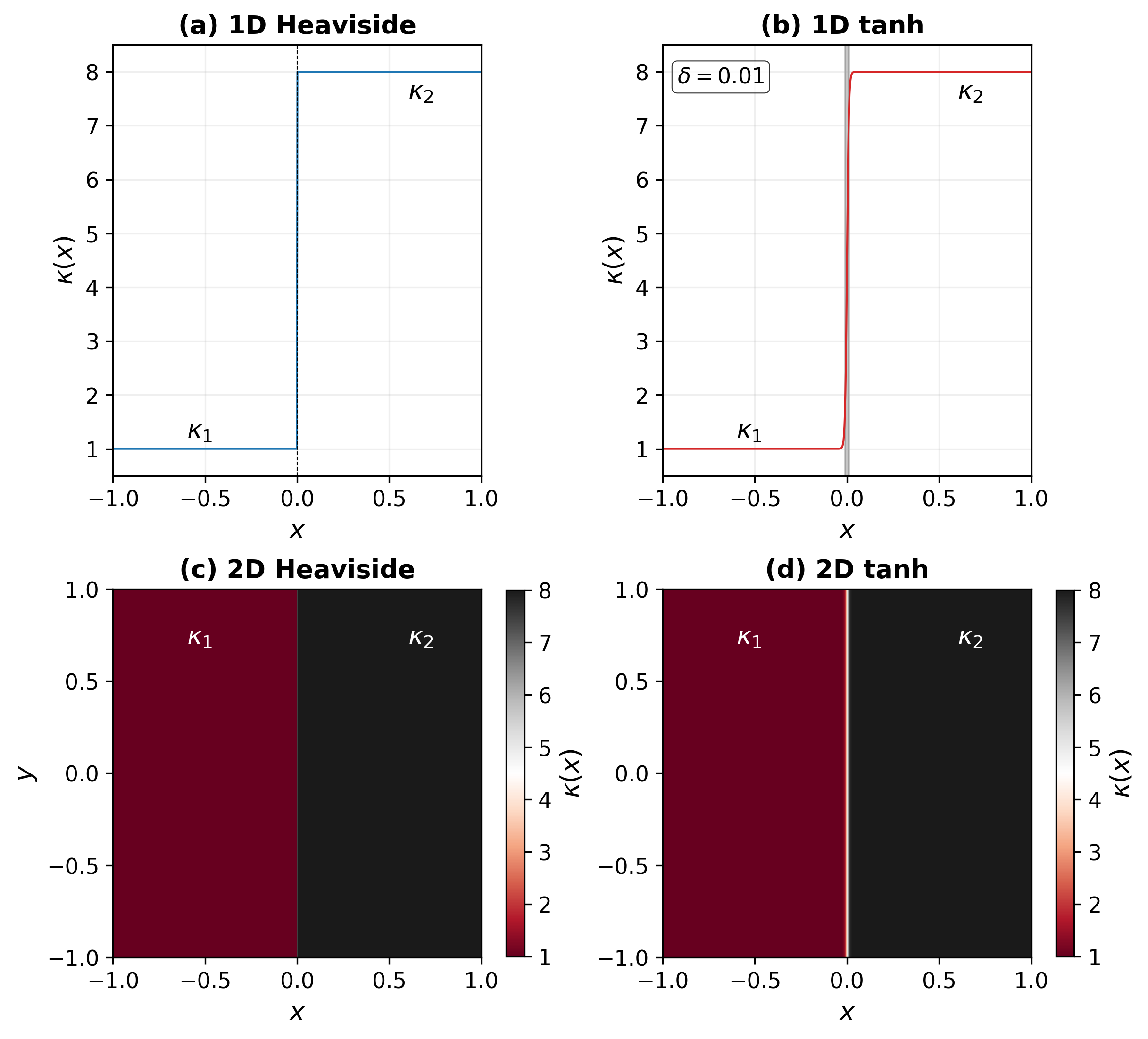}
    \caption{Material coefficient $\kappa(x)$ represented by a sharp Heaviside step 
    function (a,c) and a smooth $\tanh$ approximation with smoothing parameter $\delta=0.01$ (b,d). 
    Panels (a,b) show the 1D case, while (c,d) illustrate the corresponding 2D extensions.}
    \label{fig:Interface}
\end{figure}

By modeling $\kappa$ as a step function, as illustrated in Figure~\ref{fig:Interface}(a), the material changes suddenly at the interface $\Gamma_I$. In this case $\nabla \kappa$ behaves as a Dirac function, leading to numerical ill conditioning.

To overcome this difficulty, the problem is reformulated by introducing 
explicit interface conditions expressed with respect to the unknown function $u$: 

\begin{equation}
\begin{aligned}
[u]_{\Gamma_I} &= 0, \\
[\kappa\,\partial_n u]_{\Gamma_I} &= 0.
\end{aligned}
\label{eq:Interface}
\end{equation}

Here, $[\,\cdot\,]_{\Gamma_I}$ denotes the variation across the interface. 
The first condition enforces continuity of the function $u$ (often referred as potential function), while the 
second condition ensures conservation of the flux through $\Gamma_I$. 
Together, they replace the singular term in the strong form and provide a 
well-posed formulation for inhomogeneous materials.

In contrast, the weak form \eqref{eq:weak-local} does not require derivating $\kappa$, and therefore it remains well-defined also for non-homogeneous materials. Note that the weak form integration automatically satisfies the interface conditions of \eqref{eq:Interface}. Thus, even though $\kappa$ is discontinuous, the weak form enforces the correct interface conditions without requiring explicit penalty terms or additional constraints, and the discontinuity is handled implicitly in the integral evaluation.

\subsubsection{Smooth interface representation (tanh)}

In the strong form, to address the issue of singularities, the material 
discontinuity can be approximated by a smooth $\tanh$ function

\begin{equation}
\kappa(x_n)=\kappa_1+(\kappa_2-\kappa_1)\,T(x_n),
\qquad
T(x_n)=\tfrac12\big(1+\tanh(x_n/t)\big).
\label{eq:kappa-tanh}
\end{equation}

with the coordinate system being the same shown in Figure \ref{fig:interface_geometry} and where $t>0$ controls the thickness of the transition layer. 

Although the derivative of this function, $\tfrac{d}{dx_n}T(x_n) = \tfrac{1}{2t},\mathrm{sech}^2(x_n/t)$, is smooth and $\kappa$ remains differentiable everywhere, when \textit{t} is chosen very small to model a sharp discontinuity, numerical issues persist, and the strong form remains difficult to handle.

On the other hand, if the smooth version of $\kappa$ is used inside the weak form, the PDE remains well-defined because no derivatives of $\kappa$ appear explicitly, and the weak formulation automatically enforces the correct interface conditions. As \textit{t} $\to$ 0, the solution converges to the sharp-interface problem. 

In practice, choosing a small but finite \textit{t} regularizes the 
discontinuity in a way that preserves the underlying physics while 
avoiding numerical difficulties. As illustrated in 
Figure~\ref{fig:Interface}(b) and (d), the $\tanh$ smoothing introduces a very 
narrow but essential transition layer across the interface. 
This aspect is relevant when actually training the IVI-PINN adopting algorithms based on Automatic Differentiation (AD), such as ADAM~\cite{kingma2015adam}. As a matter of fact, using a smooth function to describe material discontinuities improves AD, since the gradients of 
$\kappa$ remain smooth and bounded, ensuring well-behaved 
sensitivities during training. It also improves the accuracy of numerical integration: in the weak form, terms such as $\int_{\Omega} \kappa(x)\,\nabla u \cdot \nabla v \, dx$  are evaluated using quadrature rules (e.g., Gauss points). If $\kappa$ is a sharp step function, the discontinuity may fall between 
quadrature points, so the jump is not captured correctly, leading to 
aliasing and sampling errors. By contrast, with a smoothed material coefficient 
$\kappa$, the integrand is continuous, so quadrature points sample 
a well-defined function and the computed integral remains accurate even 
if the points do not align with the interface.

\section{Neural network and loss function}
\label{sec:nn-loss}

From now on we refer to 2D problem, in which the element-wise representation of the 
weak form~\eqref{eq:weak-local} is considered for each element $\Omega_{k}$ 
with Neumann boundary $\Gamma_{N}^{k}$ and test functions 
$v^{(x)}_{k,n}(x)$ and $v^{(y)}_{k,m}(y)$. The corresponding bilinear form 
is obtained from the left-hand side of~\eqref{eq:weak-local}.

\begin{equation}
\begin{aligned}
U^{(k)}_{n,m} 
= & \int_{\Omega_{k}} 
\kappa(x,y)\left( \frac{\partial u(x,y)}{\partial x} \, 
       \frac{\partial v^{(x)}_{k,n}(x)}{\partial x}\, v^{(y)}_{k,m}(y)
     + \frac{\partial u(x,y)}{\partial y} \, 
       v^{(x)}_{k,n}(x)\, \frac{\partial v^{(y)}_{k,m}(y)}{\partial y}
\right) dx\,dy  \\
& - \int_{\Gamma_{N}^{k}} 
\kappa(x,y) \, g_N \, v^{(x)}_{k,n}(x)\, v^{(y)}_{k,m}(y)\, d\Gamma .
\end{aligned}
\end{equation}

The right-hand side contribution remains
\begin{equation}
F^{(k)}_{n,m} 
= \int_{\Omega_{k}}
f(x,y)\, v^{(x)}_{k,n}(x)\, v^{(y)}_{k,m}(y) \; dx\,dy.
\end{equation}

It is worth mentioning that the integrals are computed using a 
tensor–product Gauss–Lobatto–Jacobi (GLJ) quadrature rule~\cite{Doha2014_JGL_NLSE} of order 
\( N_{\rm quad} \) on each local element \( k \), with the corresponding 
quadrature weights. This procedure results in \( N_{\rm quad}^2 \) 
quadrature points per element. The use of GLJ quadrature is particularly convenient because the endpoints are included, allowing the boundary terms to be naturally captured in the weak formulation.

Thus, the residual for each test function is

\begin{equation}
\mathcal{L}_{\mathrm{V}} 
= \sum_{k=1}^{N} \frac{1}{N^{x}N^{y}}
\sum_{n=1}^{N^{x}} \sum_{m=1}^{N^{y}}
\big( R^{(k)}_{n,m} \big)^{2},
\qquad \text{in  } \Omega \cup \Gamma_N .
\end{equation}
where $R^{(k)}_{n,m} =U^{(k)}_{n,m} -F^{(k)}_{n,m}$.

For Dirichlet boundary conditions, an additional penalty term is imposed:
\begin{equation}
\mathcal{L}_{\mathrm{D}}
= \frac{1}{N_D} \sum_{N_D} 
\big( u - g_D \big)^2,
\qquad \text{in  } \Gamma_D.
\end{equation}

The total training objective for VPINNs is then
\begin{equation}
\mathcal{L}_{\mathrm{INI\_VPINN}} 
= \lambda_D\,\mathcal{L}_{D} 
+ \lambda_V\,\mathcal{L}_{\mathrm{V}},
\end{equation}
where $\lambda_D$ and $\lambda_V$ are the weights associated with the 
Dirichlet boundary loss and the weak variational loss, respectively.

\begin{figure}
    \centering
    \includegraphics[width=1\linewidth]{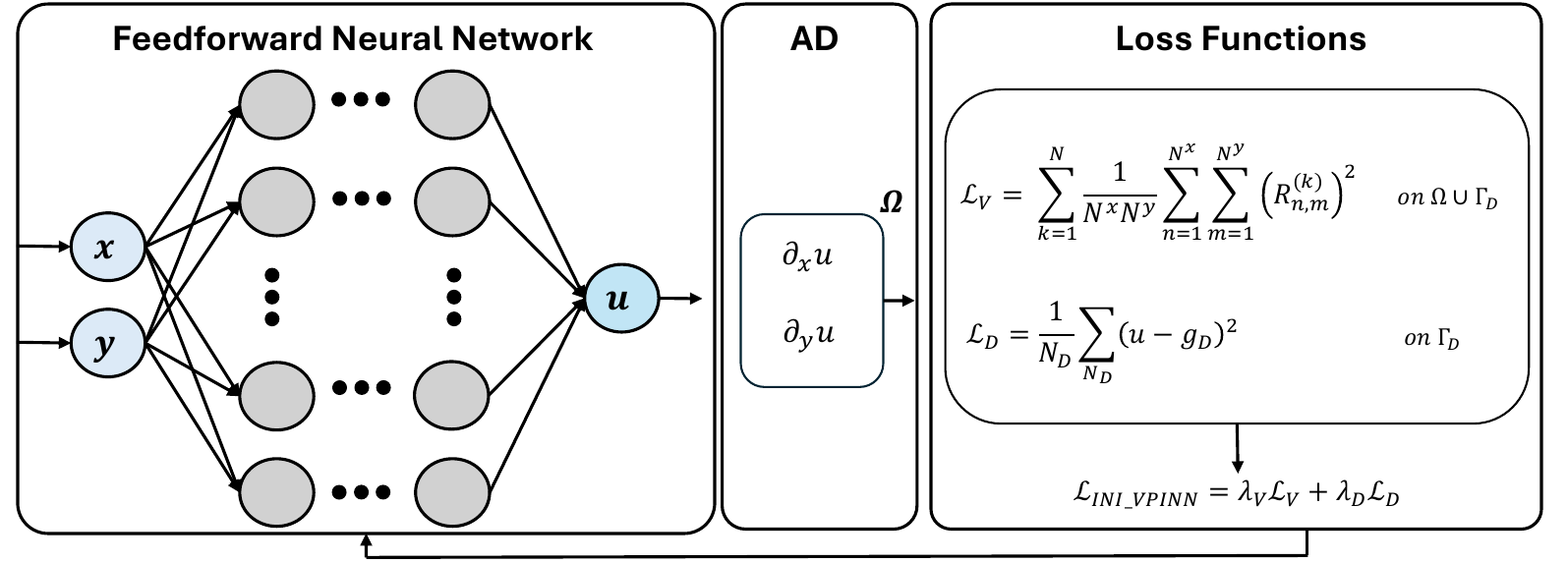}
    \caption{Schematic of the INI-VPINN framework.}
    \label{fig:NN}
\end{figure}

The fully connected neural network, sketched in Figure~\ref{fig:NN}, 
approximates $u:\Omega\to\mathbb{R}$ from the input $(x,y)$. 
The Poisson equation together with the Neumann boundary condition is enforced via a weak-form loss (mean–squared error of the element-wise variational residuals), 
while the Dirichlet conditions are imposed by a mean-squared error penalty on sampled boundary points.

\section{Results}
In this section, the performance of the INI-VPINN approach across different problem complexities, such as geometric singularity and material discontinuities, is evaluated. From a PDE perspective, the test cases can be divided into two categories: the Laplace equation and the Poisson equation with mixed boundary conditions. 

The results obtained by the INI-VPINN approach are compared with other PINN-based approaches present in the literature. The approches have been selected  to directly reflect the specific challenges addressed by the proposed
INI-VPINN formulation, namely: (i) the treatment of Neumann boundary conditions;
(ii) the handling of interface conditions in multi-material domains; (iii) the difficulty of balancing multiple loss terms in PINNs;  (iv) the presence of geometric singularities (e.g., corners), which are known to degrade the performance of standard PINN formulations.
In particular, the comparison strategy is as follows:
\begin{itemize}
\item Vanilla PINN: selected as a baseline reference, as it represents the standard formulation. It highlights known challenges related to balancing multiple loss terms (PDE,
Dirichlet, and Neumann), especially in problems with mixed boundary conditions.
\item VPINN with explicit Neumann enforcement: selected to enable
a direct comparison within a variational framework. Its weak formulation, which avoids
second-order derivatives and introduces element-wise weighting functions via domain de-
composition, facilitates assessment of the impact of explicit versus implicit enforcement of Neumann boundary conditions.
\item Conservative PINN (cPINN): selected as a representative method for multi-material
problems, where interface conditions are enforced explicitly through domain decomposition. This enables a direct comparison with the proposed implicit interface treatment in
INI-VPINN.
\end{itemize}
Therefore, due to the nature of INI-VPINN, which is designed to address multiple challenges, the comparison is performed against representative and robust PINN variants
tailored to each specific aspect. This strategy enables a structured and meaningful comparison, allowing the advantages of the proposed formulation to be clearly demonstrated.
The proposed test cases are then purposely designed to target specific challenges, and the most appropriate PINN methods, including INI-VPINN, are evaluated for each case to assess their relative performance.

In order to make a comparison with an accurate solution, all the test cases are also modeled and solved by the use of a commercial Finite Element Method (FEM) based software (i.e. Comsol Multiphysics Ver. 6.4 (c) ): quadratic elements are selected and the number of degrees of freedom (influenced by the mesh density) is selected in order to obtain constant accuracy given additional density increase.

\begin{figure} 
    \centering
    \includegraphics[width=0.7\linewidth]{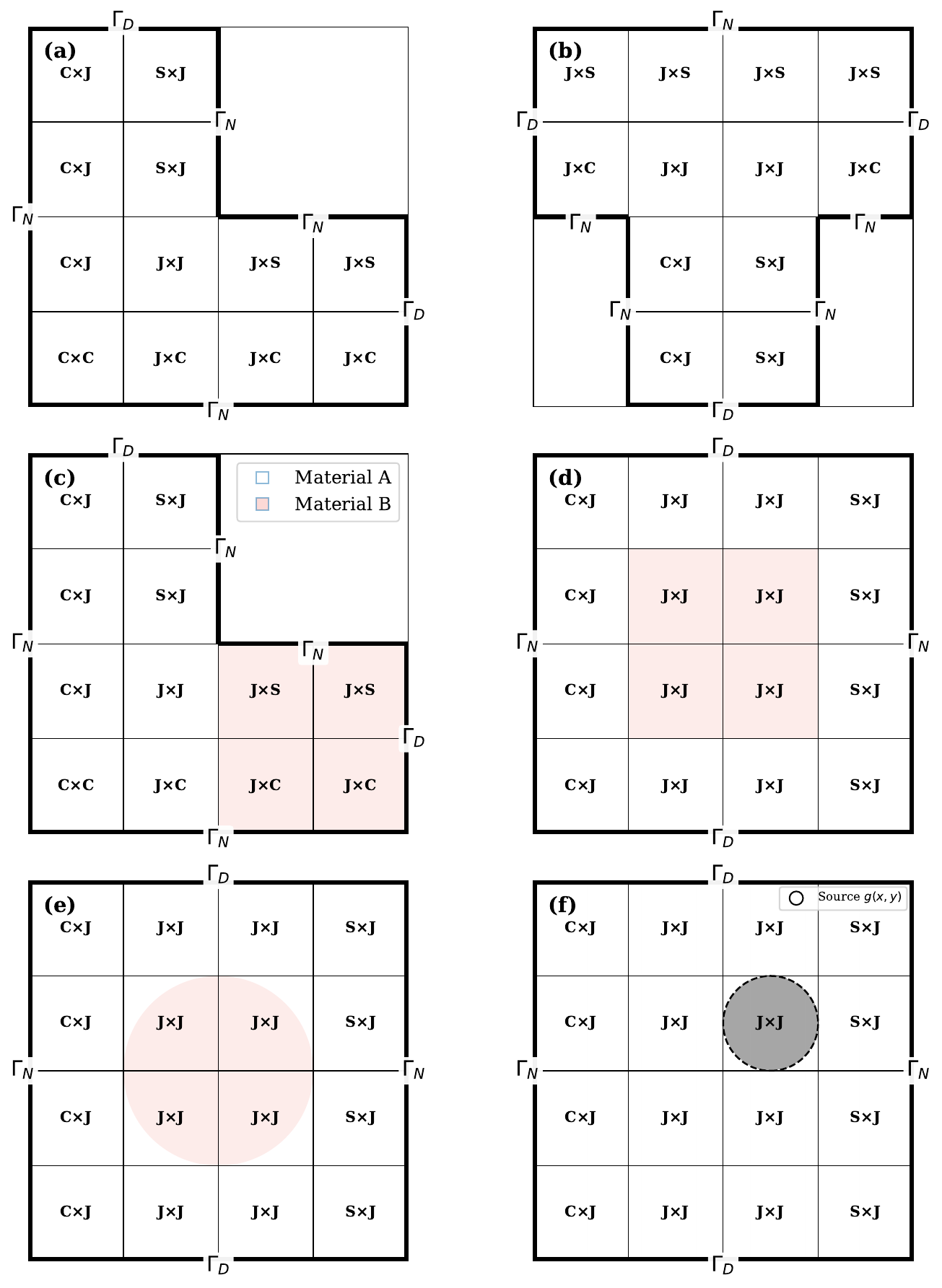} 
    \caption{Schematic representation of the geometry of the test cases and distribution of the test functions within local elements for each case.
    (a) Homogeneous L-shaped domain.
    (b) Homogeneous T-shaped domain.
    (c) Non-homogeneous L-shaped domain with two materials.
    (d–e) Non-homogeneous square domains with different material inclusions: (d) square inclusion and (e) circular inclusion.
    (f)  Square domain with Gaussian source term.}\label{fig:domains}
\end{figure}

First, two test cases with homogeneous materials but geometric singularities are examined, as illustrated in Figures~\ref{fig:domains}(a) and (b), corresponding to the L-shaped and T-shaped domains. The square circumscribing the above defined domains is characterized by a side length $L=1m$; the choice of these shapes arises from the fact that corners pose a challenge for numerical solutions, especially when they are characterized by boundary conditions of the Neumann type.

Since another strength of the INI-VPINN is its ability to handle material discontinuities, three test cases with material discontinuity are studied: (i) a L-shaped domain with heterogeneous materials (see Figure~\ref{fig:domains}(c)), (ii) a nested square configuration where the inner square has a different material from the outer region (see Figure~\ref{fig:domains}(d)), and (iii) a square domain containing a circular inclusion with different material properties (see Figure~\ref{fig:domains}(e)). The main aim of the third case is assessing IVI-PINN in scenarios where the material interface is not aligned with the test function boundaries. Finally, a Poisson equation case with mixed boundary conditions is also considered to assess the method in the presence of a source term.

The neural network employs \textit{tanh} activation functions in its hidden layers, with Xavier initialization applied to ensure stable gradient flow and mitigate vanishing or diverging values.

As for the optimization, the Adam optimizer is used with a fixed learning rate of 0.001 to minimize the overall loss, composed by the weighted sum of the variational form term and  the Dirichlet boundary condition term.

The training points used within each element correspond to the quadrature nodes of a tensor–product Gauss–Lobatto–Jacobi (GLJ) rule of order \( N_{\rm quad}=10 \), which are employed to numerically evaluate the integrals appearing in the variational (weak) formulation, resulting in \( N_{\rm quad}^2 = 100\) quadrature points per element. In addition to these points, there are \( N_{\rm Dir}=80 \) boundary points used to enforce the Dirichlet boundary condition.

The rank \( n \) of the weighting functions is selected based on the complexity of the problem. In all test cases considered here \( n = 4 \) was found to be sufficient, as higher values do not yield noticeable gains in accuracy but only add extra computational cost.

\subsection{Homogeneous L-Shaped Domain}
\label{sec:hlshape}

In this test case, we consider the Laplace equation ($f = 0$) on a L-shaped domain with mixed Neumann-Dirichlet boundary conditions. The Dirichlet boundary condition ($\Gamma_{D}$) is prescribed as $u = 1$ on the top edge and $u = 0$ on the right edge, while the remaining edges are assigned Neumann conditions ($\Gamma_{N}$: $\partial_n u = 0$), as illustrated in Figure~\ref{fig:domains}(a). The material property in this case is homogeneous, with $\kappa = 1$ assigned throughout the entire domain.

The main reason for choosing this geometry is that the corner at $(x,y)=(0,0)$ poses challenges for numerical methods, so it makes this test case an essential benchmark for evaluating the robustness of the proposed approach.

The L-shaped domain is discretized into $12$ uniform non-overlapping elements, with $N_{\rm quad} = 10$ quadrature nodes in each element.
As described in Section~\ref{sec3}, the choice of test functions in each element depends on the block location, as illustrated in Figure~\ref{fig:domains}(a). Along each direction $x$ or $y$, if the left edge corresponds to a Neumann boundary, odd cosine test functions are used; if its left edge corresponds to a Neumann boundary, odd sine test functions are adopted. Otherwise, Jacobi polynomials are the most suitable choice for interior regions or Dirichlet boundaries. 

The details of three different PINN-based models compared in this section are summarized in Table \ref{tab:LD}. The standard PINN employs a loss function composed by a (weighted) sum of the residual in the strong-form PDE and the Dirichlet and Neumann terms. The VPINN (Ex) replaces the strong-form term with its variational form. Finally, the proposed INI-VPINN further simplifies the formulation by using only the variational and Dirichlet losses, implicitly enforcing the Neumann condition within the variational framework. The number of iterations for all cases is fixed to 40,000 epochs.

\begin{table*}[h]
\centering
\small
\begin{tabular}{l l l r l r}
\toprule
\textbf{Case} & \textbf{Network} & \textbf{Loss terms} & \textbf{Iterations} \\
\midrule
PINN & $[2,5,5,5,1]$ & $\mathcal{L}_{\mathrm{S}},\; \mathcal{L}_{\mathrm{D}},\; \mathcal{L}_{\mathrm{N}}$ \,(3 terms) & 40{,}000 \\
VPINN (Ex) & $[2,5,5,5,1]$ & $\mathcal{L}_{\mathrm{V}},\; \mathcal{L}_{\mathrm{D}},\; \mathcal{L}_{\mathrm{N}}$ \,(3 terms) & 40{,}000 \\
INI-VPINN& $[2,5,5,5,1]$ & $\mathcal{L}_{\mathrm{V}},\; \mathcal{L}_{\mathrm{D}}$ \,(2 terms) & 40{,}000 \\
\bottomrule
\end{tabular}
\caption{Comparison of different models compared in the homogeneous L-shaped domain problem. The loss components correspond to $\mathcal{L}_{\mathrm{D}}$: Dirichlet boundary condition, $\mathcal{L}_{\mathrm{V}}$: Variational form, $\mathcal{L}_{\mathrm{N}}$: Neumann boundary condition, and $\mathcal{L}_{\mathrm{S}}$: Strong form.}
\label{tab:LD}
\end{table*}

\begin{figure}[H] 
    \centering
    \includegraphics[width=0.8\linewidth]{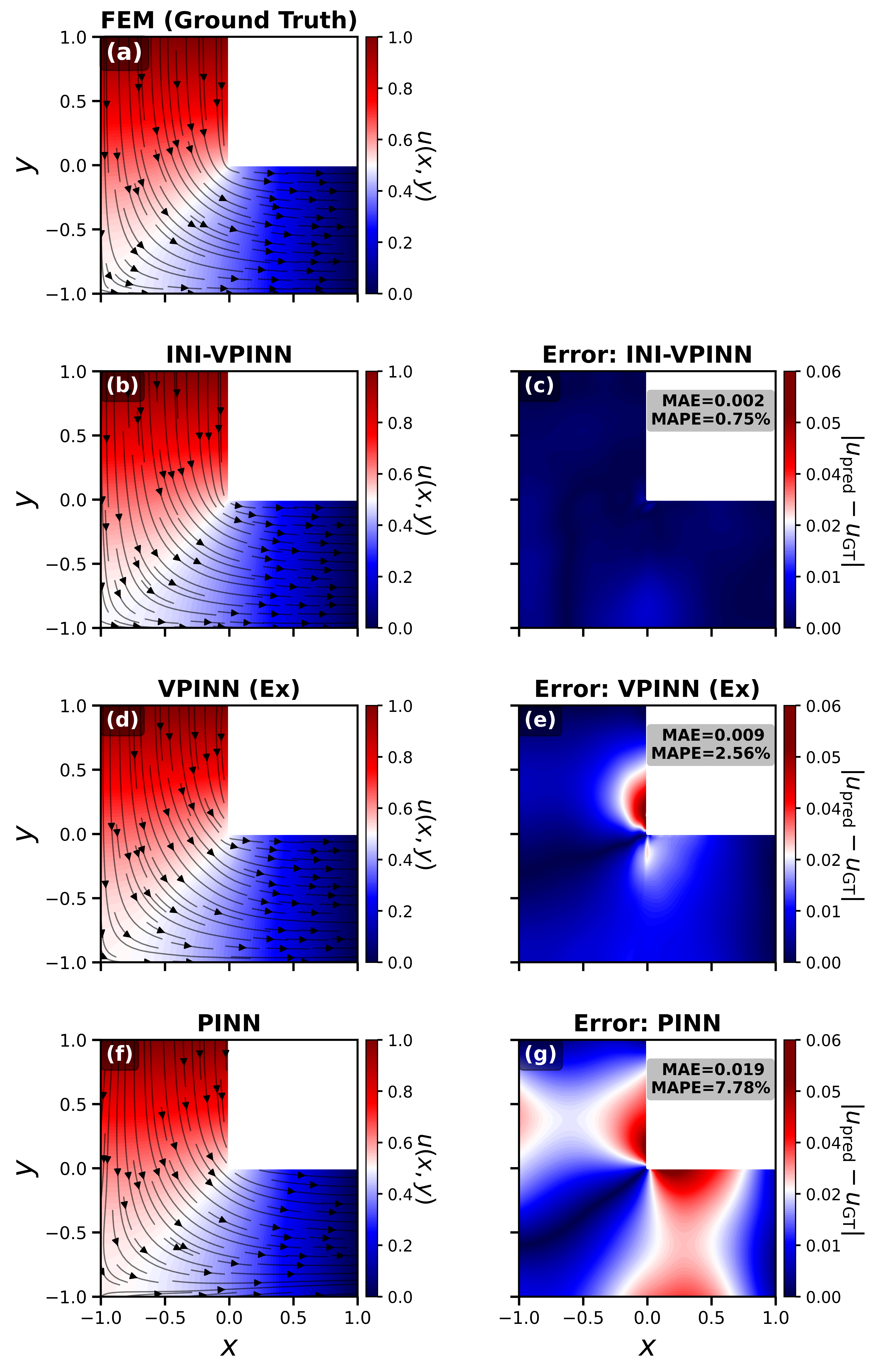} 
    \caption{
    Comparison of predicted $u(x,y)$ for the L-shaped domain with homogeneous material: 
    (a) FEM (Ground Truth), (b) INI-VPINN, (d) VPINN (Ex), with explicitly enforced Neumann boundary condition, 
    (f) PINN. Streamlines in each case illustrate the direction of the field ($-\nabla u(x,y)$). (c-e-g) Corresponding absolute error distributions $|u_{\text{Pred}} - u_{\text{GT}}|$.}
    \label{fig:LD}
\end{figure}

As shown in Figure~\ref{fig:LD} all methods reproduce the general potential distribution, but their accuracy differs significantly. INI-VPINN provides the closest agreement with FEM, achieving not only the lowest Mean Absolute Error (${\rm MAE}=0.002$, ${\rm MAPE}=0.75\%$), but also an error map that is uniformly distributed across the domain (as shown in Figure~\ref{fig:LD}c)). VPINN (Ex), in which the Neumann boundary condition is explicitly enforced, also shows good agreement, but with slightly larger localized errors near the Neumann boundaries and in the vicinity of the corner (${\rm MAE}=0.009$, ${\rm MAPE}=2.56\%$).

In contrast, the standard PINN struggles with geometries containing corners, especially in the presence of mixed boundary conditions. In this case, the loss function consists of three components—one for the PDE 
residual, one for the Dirichlet boundary condition, and one for the Neumann boundary condition. Balancing the relative contribution of these terms so that each one of them plays a significant role is the main challenge, which results in larger errors concentrated near the corner and along Neumann boundaries (${\rm MAE}=0.019$, ${\rm MAPE}=7.78\%$).

Note that while the accuracy of the solution in terms of the unknown function $u(x,y)$ is fundamental, also the derivative of $u$ has often practical importance; for instance the gradient of the potential (electric or magnetic) gives the magnitude of the field (electric or magnetic); for this reason the field lines of $(-\nabla u)$  are also shown in the results for all cases.

Field lines accuracy is especially critical under Neumann boundary conditions, where the flux must remain tangent to the boundary, and near the corner singularity, where gradients are large. As shown in the Figure~\ref{fig:LD}, INI-VPINN streamlines are well aligned with the FEM reference. On the other hand VPINN (Ex) demonstrates slight localized deviations, especially in Neumann boundaries close to the corner. Finally, PINN shows noticeable mismatches with non-tangential field lines along Neumann boundaries and evident divergence near the corner. Overall, INI-VPINN demonstrates the most robust performance, outperforming VPINN (Ex) and significantly surpassing the standard PINN in both accuracy and physical consistency.

\begin{figure}[H]
    \centering
    \includegraphics[width=0.55\linewidth]{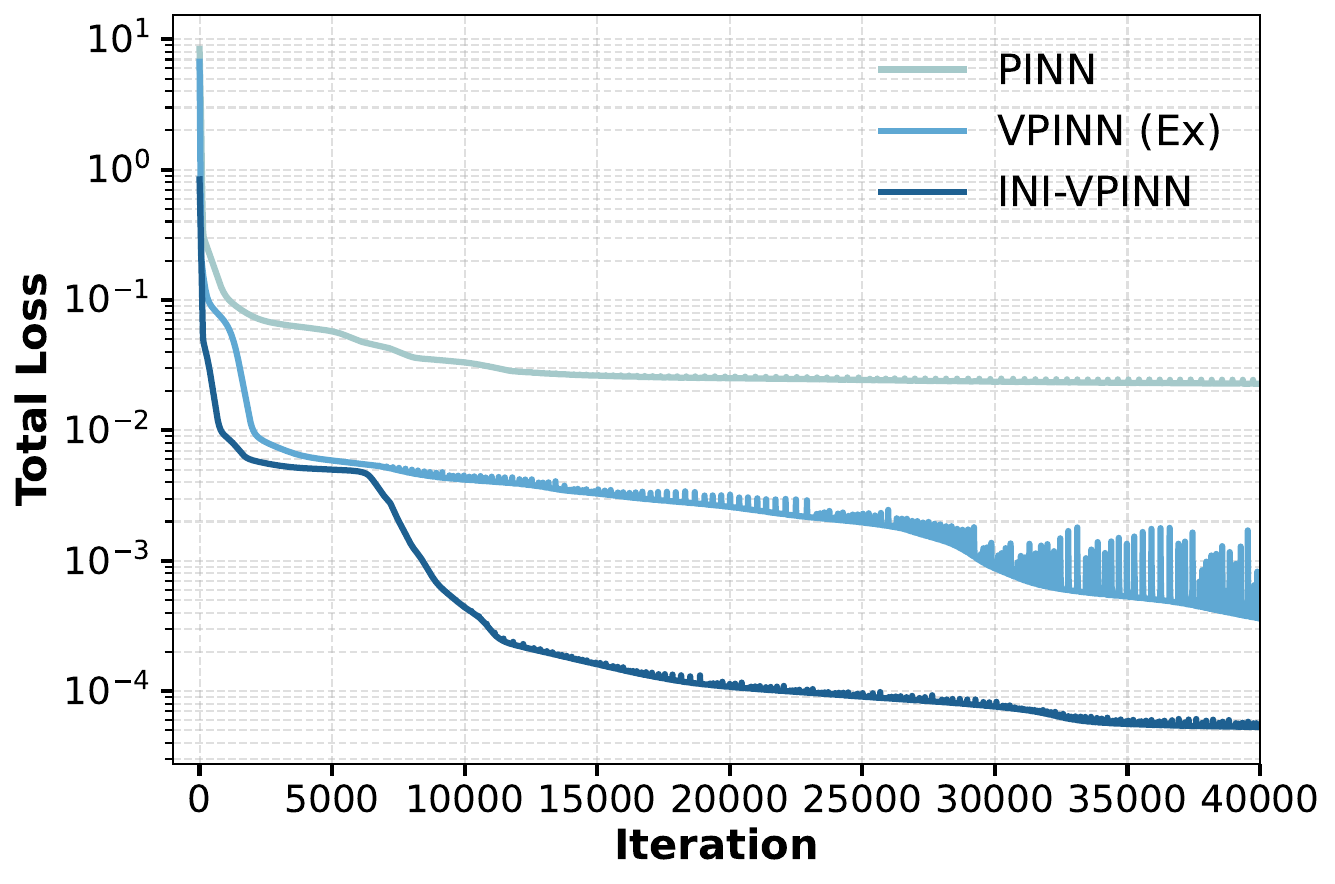}
    \caption{Training loss convergence of PINN, VPINN (Ex), and INI-VPINN models for the homogeneous L-shaped domain.}
    \label{fig:LossHistory_Comparison_1}
\end{figure}

The convergence trends shown in Figure~\ref{fig:LossHistory_Comparison_1} clearly indicate that the proposed INI-VPINN total loss curve drops much more steeply during the first few thousand iterations, which means it trains faster and more effectively among all models. In contrast, the standard PINN and VPINN (Ex) converge more slowly. All models were trained for a sufficient number of iterations to ensure convergence of their respective losses, so under these comparable conditions, the INI-VPINN clearly demonstrates faster convergence and efficiency.

\subsection{Homogeneous T-Shaped Domain}

A more challenging test case is the T-shaped geometry, characterized by mixed Neumann and Dirichlet boundary conditions and two corner singularities. As in the previous test, the Laplace equation ($f = 0$) with mixed Neumann–Dirichlet boundary conditions is considered. The Dirichlet boundary condition ($\Gamma_{D}$) is prescribed as $u = 1$ on the left, $u = 0.5$ on the right, and $u = 0$ on the bottom boundaries, while the remaining edges are subject to Neumann conditions ($\Gamma_{N}$: $\partial_n u = 0$), as illustrated in Figure~\ref{fig:domains}b).

The architecture, loss terms, and number of epochs of the two VPINN-based models used in this problem are summarized in Table \ref{tab:TD}.

\begin{table*}[h]
\centering
\small
\begin{tabular}{l l l r}
\toprule
\textbf{Case} & \textbf{Network} & \textbf{Loss terms} & \textbf{Iterations} \\
\midrule
VPINN (Ex) & $[2] + 4 \times [19] + [1]$ & $\mathcal{L}_{\mathrm{V}},\; \mathcal{L}_{\mathrm{D}},\; \mathcal{L}_{\mathrm{N}}$ \,(3 terms) & 80{,}000 \\
INI\textendash VPINN& $[2] + 4 \times [19] + [1]$ & $\mathcal{L}_{\mathrm{V}},\; \mathcal{L}_{\mathrm{D}}$ \,(2 terms) & 80{,}000 \\
\bottomrule
\end{tabular}
\caption{Comparison of different models implemented in the homogeneous T-shaped domain problem. The loss components correspond to $\mathcal{L}_{\mathrm{D}}$: Dirichlet boundary condition, $\mathcal{L}_{\mathrm{V}}$: Variational form, and $\mathcal{L}_{\mathrm{N}}$: Neumann boundary condition.}
\label{tab:TD}
\end{table*}

\begin{figure}[H] 
    \centering
    \includegraphics[width=.885\linewidth]{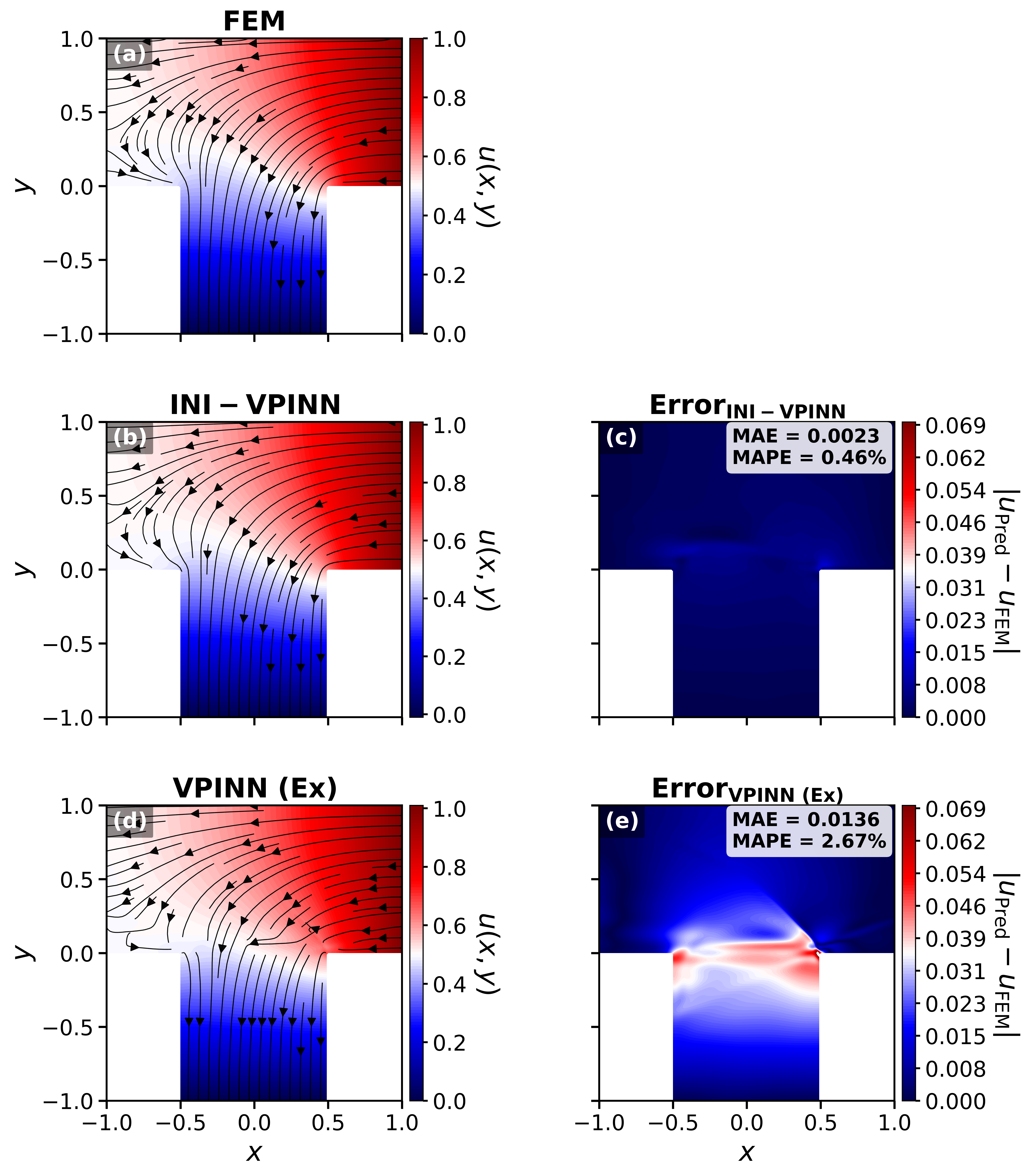} 
    \caption{Comparison of predicted $u(x,y)$ in the T-shaped domain with homogeneous material using (a) FEM (Ground Truth), (b) INI-VPINN, and (c) VPINN (Ex). The colormap represents the scalar potential, while the streamlines illustrate the flux direction $(-\nabla u)$. Error distributions relative to FEM are shown in (c–e).
}\label{fig:TD}
\end{figure}

As shown in Figure~\ref{fig:TD}, the VPINN (Ex) provides a reasonable approximation of the solution over the whole domain (${\rm MAPE}=2.67\%$), but the localized errors near the corners and along the Neumann boundaries are significant. By contrast, INI-VPINN achieves significantly improved accuracy, yielding both the lowest mean error (${\rm MAE}=0.0023$, ${\rm MAPE}=0.46\%$) and an error map that is uniformly distributed across the domain. 

The INI-VPINN field lines, more sensitive to slight deviations in the potential, remain well aligned with the FEM reference, especially near corners. VPINN (Ex), on the other hand, shows visible deviations near reentrant corners, where flux lines fail to remain tangential to Neumann boundaries. Overall, INI-VPINN again demonstrates superior performance, providing the most consistent agreement with FEM in both solution accuracy and flux representation.

\begin{figure}
    \centering
    \includegraphics[width=0.55\linewidth]{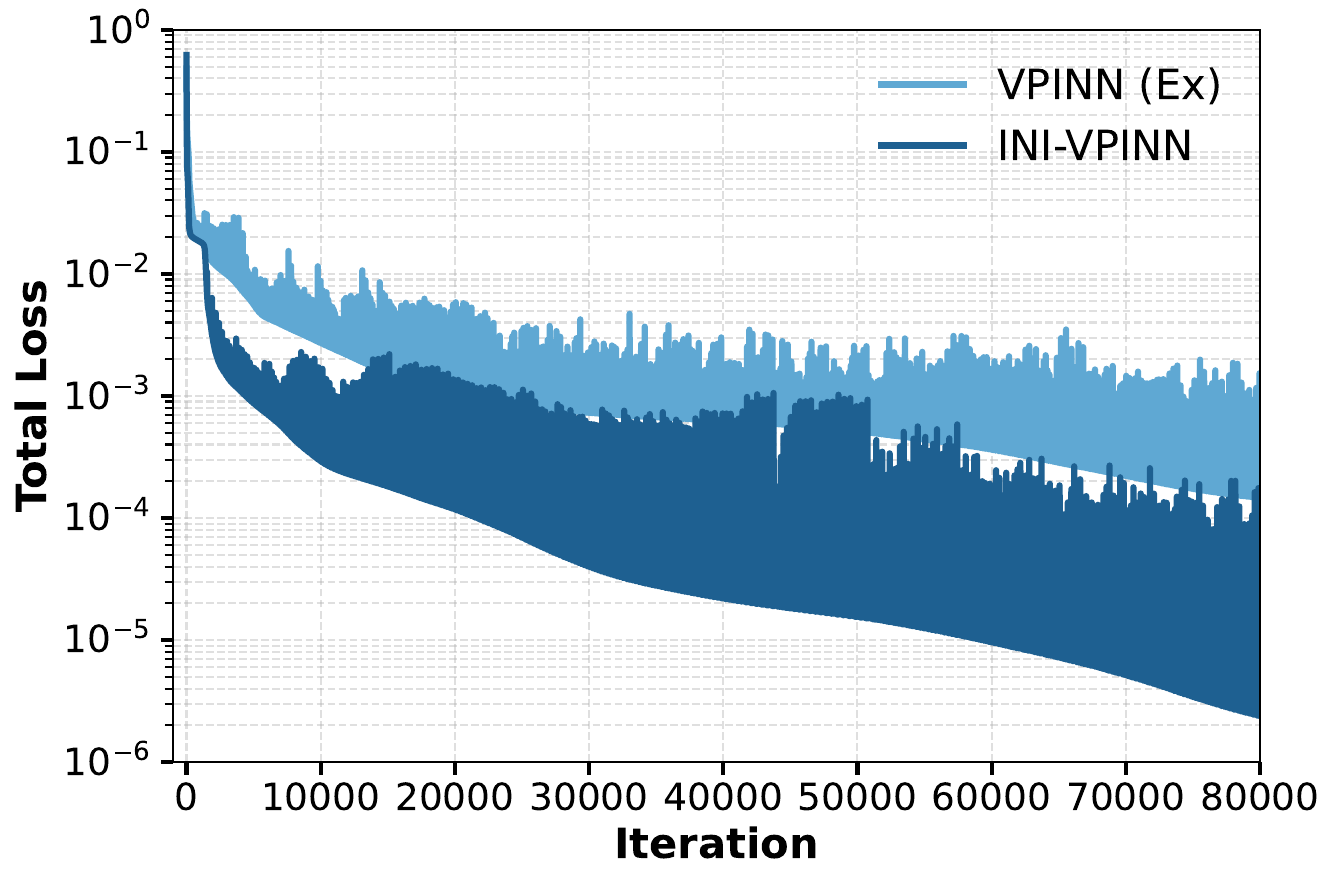}
    \caption{Training loss convergence comparison between VPINN (Ex) and INI-VPINN for the homogeneous T-shaped domain.}
    \label{fig:LossHistory_Comparison_2}
\end{figure}

The convergence behavior of the two compared models, presented in Figure~\ref{fig:LossHistory_Comparison_2}, clearly reveals the better performance of the proposed INI-VPINN compared to the explicit VPINN. Although both models show a gradual decrease in total loss, the INI-VPINN converges faster, eventually achieving a final figure two orders of magnitude smaller than the VPINN.

\subsection{Non-Homogeneous L-Shaped Domain}

In this test case, the Laplace equation ($f=0$) on a non-homogeneous L-shaped domain with a material discontinuity at $x=0$ is considered. The presence of two materials introduces additional challenges for learning-based methods. The material property for the left half of the domain ($x<0$) is set to $\kappa=1$, while for the right half ($x>0$) it is set to $\kappa=8$. In addition, mixed Neumann–Dirichlet boundary conditions are imposed as described in  section~\ref{sec:hlshape}. 

\begin{figure}[H] 
    \centering
    \includegraphics[width=.9\linewidth]{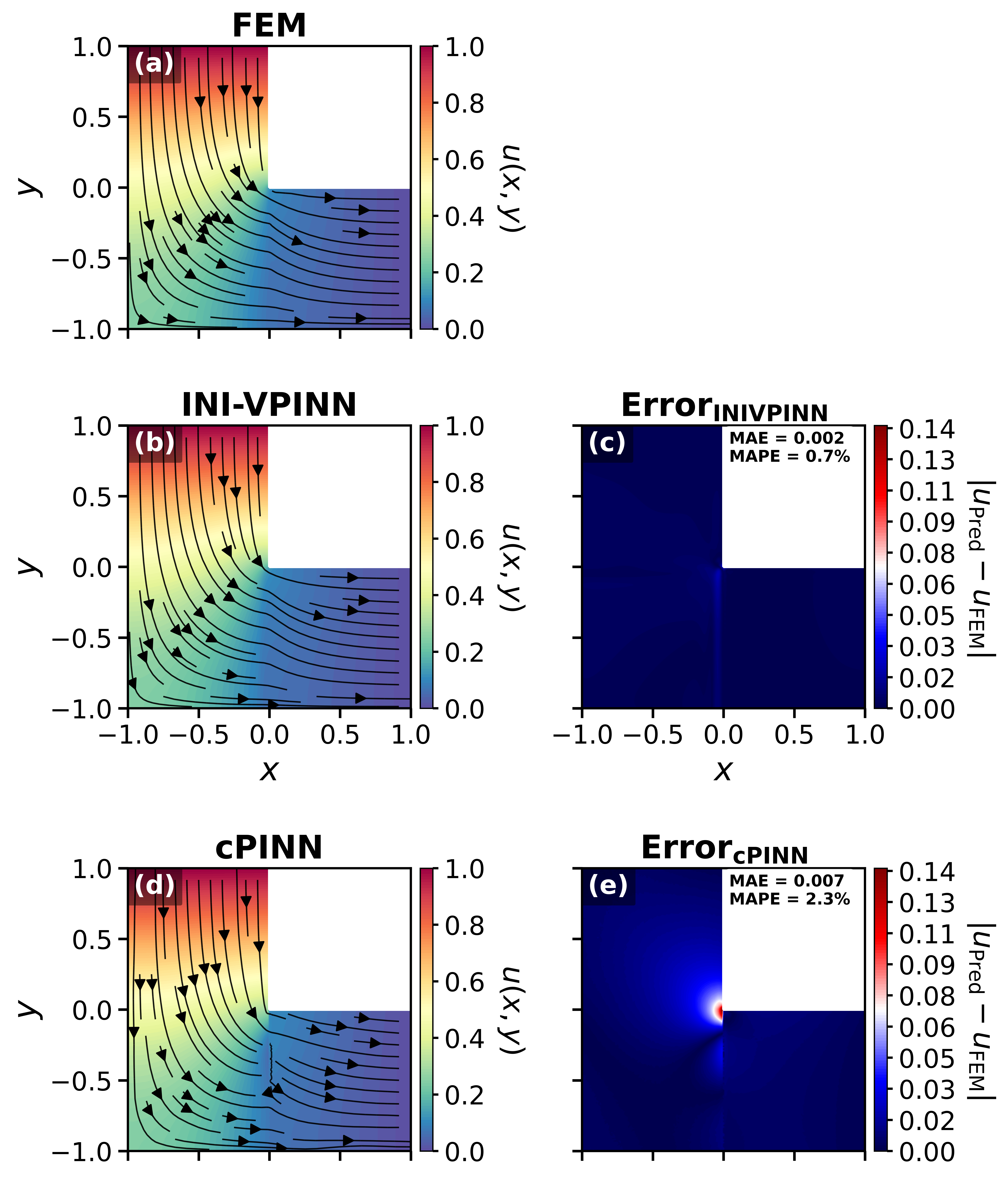} 
    \caption{Comparison of predicted $u(x,y)$ in the non-homogeneous L-shaped domain with two materials, 
            where $\kappa = 1$ for $x < 0$ and $\kappa = 8$ for $x > 0$. The solutions obtained with (a) FEM (Ground Truth), (b) INI-VPINN, and (d) cPINN (conservative physics-informed neural network) are shown. The colormap represents the scalar potential, while the streamlines illustrate the flux direction $(-\nabla u)$. Error distributions relative to FEM are 
            depicted in (c–e).}\label{fig:NHL}
\end{figure}

To evaluate the ability to handle material discontinuities in the presence of geometric singularities, we compare IVI-VPINN with the conservative physics-informed neural network (cPINN)~\cite{jagtap2020conservative}  as described in Figure~\ref{fig:NHL}. The cPINN framework is specifically designed for multi-subdomain problems, where each subdomain is assigned its own neural network and continuity across interfaces is enforced through additional loss terms.

As illustrated in Table~\ref{tab:Arc_NHLD}, the cPINN model employs two separate networks to handle the two materials, with each network consisting of three hidden layers containing 32 neurons. In contrast, the proposed INI-VPINN uses a single network with four hidden layers of 16 neurons each. Moreover, the cPINN requires balancing four different loss terms with appropriate weighting, whereas the INI-VPINN formulation involves only two loss terms, simplifying the training process.

\begin{table*}[h]
\centering
\small
\begin{tabular}{l l l r}
\toprule
\textbf{Case} & \textbf{Network} & \textbf{Loss terms} & \textbf{Iterations} \\
\midrule
INI-VPINN  & $[2] + 4 \times [16] + [1]$ & $\mathcal{L}_{\mathrm{V}},\; \mathcal{L}_{\mathrm{D}}$ & 70{,}000 \\
cPINN & 
\begin{tabular}[c]{@{}l@{}}$[2] + 3 \times [32] + [1]$ \\ $[2] + 3 \times [32] + [1]$\end{tabular} &
$\mathcal{L}_{\mathrm{S}},\; \mathcal{L}_{\mathrm{D}},\; \mathcal{L}_{\mathrm{N}},\; \mathcal{L}_{\mathrm{Int}}$ & 70{,}000 \\
\bottomrule
\end{tabular}
\caption{Comparison of different models implemented in the non-homogeneous L-shaped domain problem. The loss components correspond to $\mathcal{L}_{\mathrm{D}}$: Dirichlet boundary condition, $\mathcal{L}_{\mathrm{V}}$: Variational form, $\mathcal{L}_{\mathrm{N}}$: Neumann boundary condition, $\mathcal{L}_{\mathrm{S}}$: Strong form, and $\mathcal{L}_{\mathrm{Int}}$: Interface condition.}
\label{tab:Arc_NHLD}
\end{table*}

This decomposition allows material discontinuities to be handled more explicitly and provides a natural baseline for the non-homogeneous case, yielding reasonable errors (${\rm MAE}=0.007$ and ${\rm MAPE}=2.30\%$). On the other hand, it still struggles near the material interface and the concave corner, where enforcing flux continuity is most challenging. As a result, the localized error concentrations, as shown in Figure~\ref{fig:NHL}c), are significant.

By contrast, INI-VPINN incorporates the interface condition implicitly through its variational formulation and integrated test functions. As a result, INI-VPINN achieves substantially higher accuracy, yielding the lowest mean error (${\rm MAE}=0.002$, ${\rm MAPE}=0.70\%$) and an error map that is uniformly distributed across the domain (Figure~\ref{fig:NHL}e). The field lines further emphasize this difference: cPINN field lines deviate near the interface and lose tangential alignment along Neumann boundaries while INI-VPINN streamlines remain consistent with the FEM reference. It can be concluded that IVI-VPINN preserves both flux continuity across the interface and accuracy near the singular corner. 

Overall, although cPINN provides a structured way to address non-homogeneous domains through element decomposition, its performance is limited by the difficulty of balancing the residual, boundary, and interface losses. INI-VPINN, by embedding interface conditions into the weak form, demonstrates superior robustness and accuracy in handling this challenging two-material problem. It is also worth noting that INI-VPINN employs a single neural network, in contrast to cPINN which requires one network per subdomain, thereby offering lower computational cost for complex problems involving many materials.

\subsection{Non-Homogeneous Domain with Square Material Discontinuity}
\label{sec:S}

In this test case, we consider the Laplace equation ($f=0$) on a non-homogeneous nested square domain, where the material property is set to $\kappa = 8$ inside the square inclusion and $\kappa = 1$ outside. This configuration introduces a sharp material interface with four inner corners, making it a challenging benchmark. In particular, the continuity of both the potential and the normal flux across the material boundary must be preserved. The geometric singularities at the corners amplify the difficulty of the problem. The Laplace equation ($f = 0$) with mixed Neumann–Dirichlet boundary conditions is considered in this case. The Dirichlet boundary condition ($\Gamma_{D}$) is prescribed as $u = 1$ on the top and $u = 0$ on the bottom boundaries, while on the remaining edges  Neumann conditions ($\Gamma_{N}$: $\partial_n u = 0$) are imposed, as illustrated in Figure~\ref{fig:domains}d).

\begin{figure}[H] 
    \centering
    \includegraphics[width=.8\linewidth]{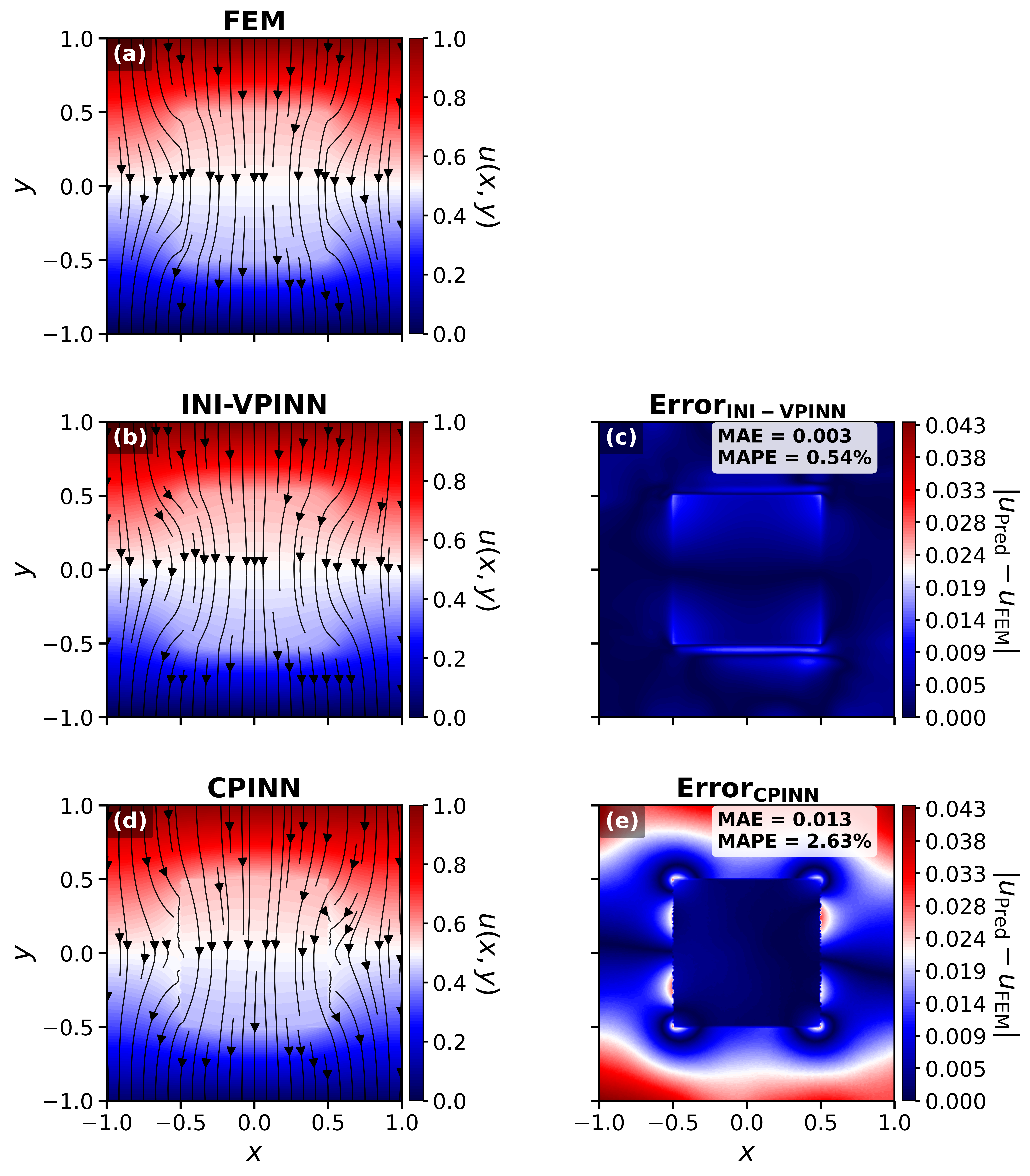}
    \caption{Comparison of predicted $u(x,y)$ in the non-homogeneous nested square domain, 
    with $\kappa = 8$ inside the square inclusion and $\kappa = 1$ outside. 
    The solutions obtained with 
    (a) FEM (Ground Truth), 
    (b) INI-VPINN, and 
    (d) cPINN (conservative physics-informed neural network) are shown. The colormap represents the scalar potential, 
    while the streamlines illustrate the field direction $(-\nabla u)$. Error distributions relative to FEM 
    are depicted in (c–e).}\label{fig:S}
\end{figure}

\begin{table*}[t]
\centering
\small
\begin{tabular}{l l l r}
\toprule
\textbf{Case} & \textbf{Network} & \textbf{Loss terms} & \textbf{Iterations} \\
\midrule
INI-VPINN  & $[2] + 5 \times [20] + [1]$ & $\mathcal{L}_{\mathrm{V}},\; \mathcal{L}_{\mathrm{D}}$ & 70{,}000 \\
cPINN & 
\begin{tabular}[c]{@{}l@{}}$[2] + 3 \times [32] + [1]$ \\ $[2] + 3 \times [32] + [1]$\end{tabular} &
$\mathcal{L}_{\mathrm{S}},\; \mathcal{L}_{\mathrm{D}},\; \mathcal{L}_{\mathrm{N}},\; \mathcal{L}_{\mathrm{Int}}$ & 70{,}000 \\
\bottomrule
\end{tabular}
\caption{Comparison of different models implemented in the non-homogeneous domain with square and circular material Discontinuity problems. The loss components correspond to $\mathcal{L}_{\mathrm{D}}$: Dirichlet boundary condition, $\mathcal{L}_{\mathrm{V}}$: Variational form, $\mathcal{L}_{\mathrm{N}}$: Neumann boundary condition, $\mathcal{L}_{\mathrm{S}}$: Strong form, and $\mathcal{L}_{\mathrm{Int}}$: Interface condition.}
\label{tab:Arc_NHSD}
\end{table*}

As shown in Figure~\ref{fig:S}, while the cPINN framework captures the general potential distribution (${\rm MAE}=0.013$, ${\rm MAPE}=2.63\%$), it produces noticeable errors around the inner corners and along the interface, indicating significant difficulties in maintaining continuity across the interfaces.

By contrast, INI-VPINN incorporates interface conditions implicitly within its variational formulation, 
allowing the method to represent both the solution and the flux more accurately across the discontinuity. 
As a result, INI-VPINN achieves a substantially lower error (${\rm MAE}=0.003$, 
${\rm MAPE}=0.54\%$) and a more uniform error distribution across the domain, as shown in 
Figure~\ref{fig:S}c). The field lines further emphasize this improvement: while cPINN flux lines deviate 
and fail to remain smooth across the material boundary, INI-VPINN streamlines remain consistent with the 
FEM reference, preserving tangential alignment along Neumann boundaries and continuity across the 
interface. 

\begin{figure}
    \centering
    \includegraphics[width=1\linewidth]{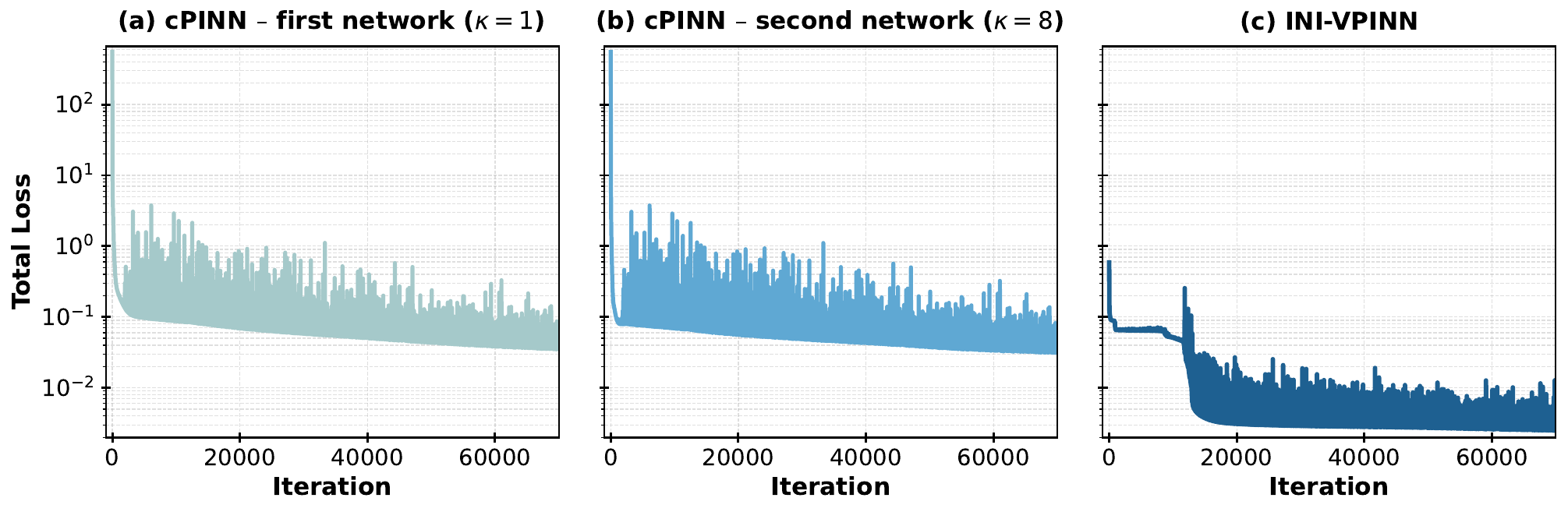}
    \caption{Convergence histories of (a) cPINN – first network ($\kappa = 1$), (b) cPINN – second network ($\kappa = 8$), and (c) the proposed INI-VPINN.}
    \label{fig:LossHistory_Comparison_3}
\end{figure}

Figure ~\ref{fig:LossHistory_Comparison_3} compares how the cPINN and INI-VPINN models converge during training. The cPINN uses two separate subnetworks—one for each material (with $\kappa = 1$ and $\kappa = 8$)—and each one works on minimizing its own loss function. However, balancing multiple loss components and maintaining interface continuity constraints remains challenging for this approach. By comparison, the proposed INI-VPINN approach converges much more smoothly and steadily, reaching lower error levels in fewer steps. All models were trained until their respective losses had stabilized, ensuring a fair comparison of their convergence characteristics.

\subsection{Non-Homogeneous Domain with Circular Material Discontinuity}

In this test case, we consider the Laplace equation ($f = 0$) on a square domain containing a circular subdomain with a different material property ($\kappa = 8$). This configuration introduces a curved interface separating two materials  making it a challenging benchmark. This is due to the lack of alignment between the subdomains interface and the elements' weighting functions at the material discontinuity interface.

The mixed Dirichlet–Neumann boundary conditions are prescribed on the outer boundaries: The Dirichlet boundary condition ($\Gamma_{D}$) on the top ($u = 1$) and on the bottom ($u = 0$) boundaries, and Neumann conditions ($\Gamma_{N}$: $\partial_n u = 0$) on the lateral boundaries, as illustrated in Figure~\ref{fig:domains}(e).

\begin{figure}[H] 
    \centering
    \includegraphics[width=.9\linewidth]{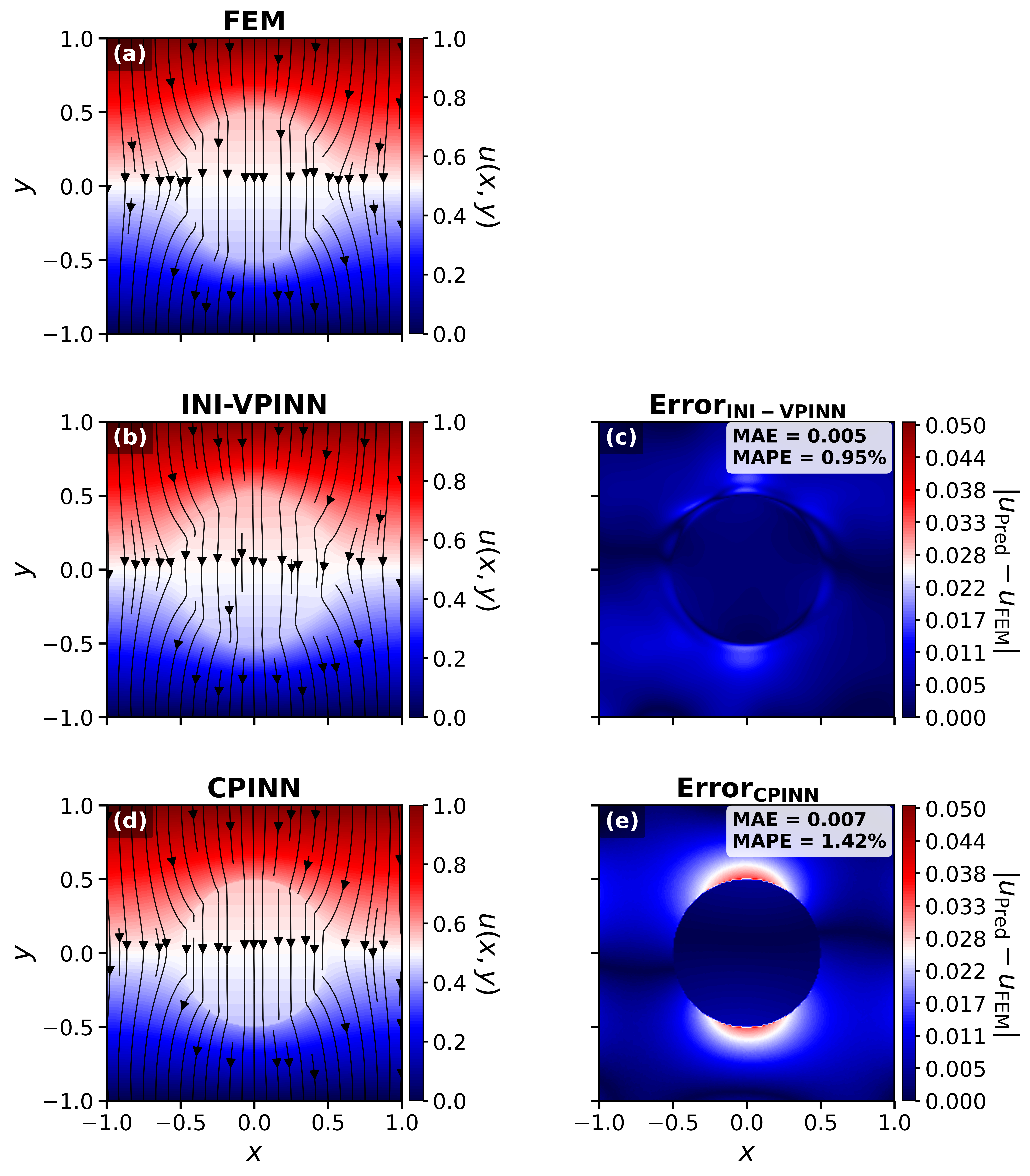}
   \caption{Comparison of the predicted potential $u(x,y)$ in the non-homogeneous square domain with circular inclusion, where $\kappa = 8$ inside the inclusion and $\kappa = 1$ outside. The solutions obtained using (a) FEM (ground truth), (b) INI-VPINN, and (d) CPINN (conservative physics-informed neural network) are presented. The colormap represents the scalar potential, while streamlines indicate the field direction ($-\nabla u$). Panels (c)--(e) show the corresponding absolute error distributions with respect to FEM.}\label{fig:O}
\end{figure}

As shown in Figure~\ref{fig:O}, the cPINN approach successfully reproduces the potential distribution in both the inner (circular) and outer subdomains, achieving ${\rm MAE} = 0.007 $ and ${\rm MAPE} = 1.42\% $. However, its main limitation lies in handling the material interface, where the largest errors occur. By contrast, the INI-VPINN formulation embeds the interface conditions implicitly within its weak-form framework,  allowing for accurate enforcement of both potential and flux continuity without explicit interface sampling. Consequently, INI-VPINN achieves superior accuracy (${\rm MAE} = 0.005$, ${\rm MAPE} = 0.95\%$) and produces smoother error distributions across the domain.

\subsection{Square domain with source term}

In this section, we evaluate the accuracy of the proposed approach in solving a Poisson problem on a square domain containing a compact source localized within the domain. The mixed boundary conditions are the same as those described in Section~\ref{sec:S}.

The source term $f(x,y)$ models a piecewise-constant source supported on a circular region centered at $(c_x,c_y) = (0.25,\,0.25)$ with radius $r=0.25$, as illustrated in Figure~\ref{fig:domains}(f):

\begin{equation}
\label{eq:source-def}
f(x,y) = A\,G(x,y), 
\qquad 
G = 
\begin{cases}
1, & (x-c_x)^2+(y-c_y)^2 \le r^2,\\[2pt]
0, & \text{otherwise},
\end{cases}
\end{equation}
where $A=-22.5$ sets the source strength.

\begin{table*}[b]
\centering
\small
\begin{tabular}{l l l r}
\toprule
\textbf{Case} & \textbf{Network} & \textbf{Loss terms} & \textbf{Iterations} \\
\midrule
INI-VPINN  & $[2] + 5 \times [20] + [1]$ & $\mathcal{L}_{\mathrm{V}},\; \mathcal{L}_{\mathrm{D}}$ & 20{,}000 \\
\bottomrule
\end{tabular}
\caption{Summary of INI-VPINN architecture and training setup in the square domain with the source term problem. The loss components correspond to $\mathcal{L}_{\mathrm{D}}$: Dirichlet boundary condition, and $\mathcal{L}_{\mathrm{V}}$: Variational form.}
\label{tab:Poison}
\end{table*}

\begin{figure}[H] 
    \centering
    \includegraphics[width=0.9\linewidth]{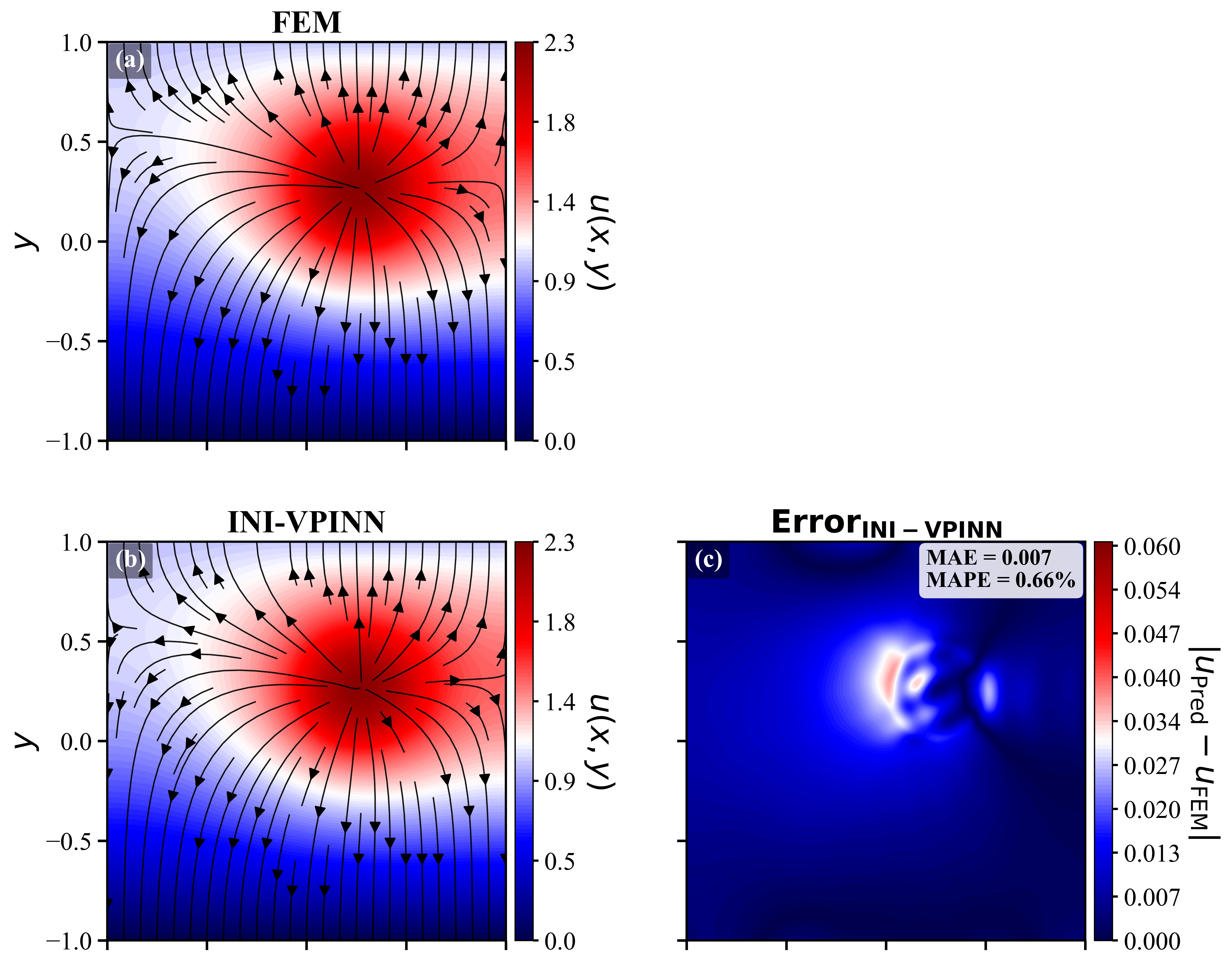}
   \caption{Comparison between the FEM reference and the proposed INI-VPINN solution for the Poisson problem with a localized source.(a) FEM reference solution (b) INI-VPINN prediction (c) Absolute error distribution.}\label{fig:source}
\end{figure}

As shown in Figures~\ref{fig:source}(a)--(b), the INI-VPINN is able to successfully predict the potential field distribution $u(x,y)$ and the corresponding streamlines of $-\nabla u$, closely matching the FEM results. The error map shows that the differences between the two solutions are very small and mainly localized near the source region, with a mean absolute error (MAE) of 0.007 and a mean absolute percentage error (MAPE) of 0.66\%.

\section{Conclusions}
This study presents the INI-VPINN, a novel weak-form Physics-Informed Neural Network approach for the resolution of Partial Differential Equations. The main feature of the approach is a specific choice of the weighting functions, leading to a native  enforcement of interface conditions and of Neumann boundary conditions thanks to its variational formulation.
The numerical experiments on Poisson and Laplace problems with sharp interfaces and geometrically complex domains confirm the accuracy and robustness of the proposed approach. Compared to several existing weak and strong form PINN variants, the INI-VPINN exhibits higher fidelity in reproducing interface behavior and geometrical singularities. The comparison between the above mentioned approaches is carried out considering a ground truth solution obtained by a classical Finite Element Method software.

It is worth mentioning here that
classical approaches, e.g., FEM or BEM, were proposed decades ago, and have nowadays achieved highly optimized performance in terms of both computational time and accuracy. Conversely, physics-guided deep learning approaches are still in their exploratory stage, and their formulations and software implementations are far from being comparable to commercial software.
Therefore, we fully agree that classical numerical methods remain the most efficient and
reliable choice, until now, for solving PDEs when the problem is fully specified, and a
high-quality mesh can be constructed. On the other hand, the goal of this work is not to
compete with or replace such well-established approaches. Instead, the aim of this contribution is to propose a new PINN-based framework (INI-VPINN) to address scenarios
where existing PINN formulations face practical limitations. In particular, the proposed
framework is advantageous in the following contexts:
\begin{enumerate}
\item \textit{ Mesh-free formulation}: Unlike FEM, INI-VPINN does not require mesh generation,
which can be challenging for complex geometries, singular domains, or problems with
evolving interfaces.
\item  \textit{Fast inference after training}: Once trained, the model acts as a surrogate solution that can be evaluated at arbitrary points in the domain with negligible computational cost. While classical methods such as FEM rely on interpolation to evaluate the solution at off-mesh locations, neural networks provide a continuous, mesh-independent representation, which can be advantageous in applications that require frequent, flexible sampling of the solution.
\item  \textit{Flexibility for inverse and data-driven problems}: Classical methods such as FEM are primarily tailored for forward problems and often require substantial reformulation for inverse settings. In contrast, PINN-based approaches naturally combine physics and data in a unified framework. Although not the focus of this work, INI-VPINN has the potential to be extended to problems involving unknown parameters or limited data.
\end{enumerate}
We emphasize that FEM solutions are used as reference benchmarks throughout the
paper, and the results demonstrate that INI-VPINN achieves comparable accuracy. In
addition, IVI-VPINN shows increased flexibility and improved accuracy and convergence in complex scenarios compared to previous PINN variants. Therefore, INI-VPINN should be regarded as a complementary approach, offering improved effectiveness compared to existing PINN variants, rather than as a replacement for traditional numerical methods.

Future developments will focus on extending the INI-VPINN to transient and nonlinear PDEs, exploring adaptive strategies for selecting weighting functions and quadrature schemes.
In summary, the INI-VPINN provides a mathematically grounded, physically consistent, and computationally robust framework for solving PDEs with complex geometries and boundary conditions. By unifying interface treatment and boundary enforcement within a single variational formulation, it represents a significant step toward the practical deployment of Physics-Informed Neural Networks in multi-physics and engineering applications.

\bibliographystyle{unsrtnat}
\bibliography{references}  






\end{document}